\documentclass[english, 12pt]{article}
\usepackage{a4wide}
\usepackage{lmodern}
\usepackage[T1]{fontenc}
\usepackage[latin9]{inputenc}
\usepackage{float}
\usepackage{bm}
\usepackage{amsmath}
\usepackage{amssymb}
\usepackage{graphicx}

\makeatletter

\providecommand{\tabularnewline}{\\}

\usepackage{xcolor}
\usepackage{mdframed}
\usepackage{marginnote}
\usepackage{todonotes}

\usepackage{subfig}
\makeatother

\usepackage{babel}
\usepackage[affil-it]{authblk}

\begin{document}

\newcommand{\barann}{\bar{A}_{nn}}
\newcommand{\barh}{\bar{h}}
\newcommand{\bfn}{\boldsymbol n}
\newcommand{\bfI}{\boldsymbol I}
\newcommand{\bfP}{\boldsymbol P}
\newcommand{\bfQ}{\boldsymbol Q}
\newcommand{\bfx}{\boldsymbol x}
\newcommand{\bfa}{\boldsymbol a}
\newcommand{\bfb}{\boldsymbol b}
\newcommand{\bfp}{\boldsymbol p}
\newcommand{\bfy}{\boldsymbol y}
\newcommand{\vel}{\boldsymbol{\beta}}
\newcommand{\conc}{u}
\newcommand{\intf}{\Sigma}
\newcommand{\bfv}{\boldsymbol v}
\newcommand{\bfu}{\boldsymbol u}
\newcommand{\bfw}{\boldsymbol w}
\newcommand{\bft}{\boldsymbol t}
\newcommand{\bfB}{\boldsymbol B}
\newcommand{\bfbeta}{\boldsymbol \beta}
\newcommand{\bfrho}{\boldsymbol \rho}
\newcommand{\bfkappa}{\boldsymbol \kappa}
\newcommand{\bfzero}{\boldsymbol 0}
\newcommand{\bfX}{\boldsymbol X}
\newcommand{\bfV}{\boldsymbol V}
\newcommand{\bfpi}{{\boldsymbol \pi}}

\newcommand{\uh}{u_h}
\newcommand{\Vsp}{{\mathcal{V}}}
\newcommand{\Vcp}{{\mathcal{V}}_c^p}
\newcommand{\Vdo}{{\mathcal{V}}_d^0}
\newcommand{\Vdp}{{\mathcal{V}}_d^p}
\newcommand{\Wspace}{{\mathcal{V}_h}
}\newcommand{\Wspaceloc}{\mathcal{Q}_{p_K}(K)}
\newcommand{\mcP}{\mathcal{P}}
\newcommand{\mcQ}{\mathcal{Q}}
\newcommand{\mcV}{\mathcal{V}}
\newcommand{\mcK}{\mathcal{K}}
\newcommand{\mcE}{\mathcal{E}}
\newcommand{\mcN}{\mathcal{N}}
\newcommand{\mcF}{\mathcal{F}}
\newcommand{\mcA}{\mathcal{A}}
\newcommand{\mcT}{\mathcal{T}}
\newcommand{\mcR}{\mathcal{R}}
\newcommand{\mcH}{\mathcal{H}}
\newcommand{\cspx}{{\mathcal{C}\mathcal{P}}_X}
\newcommand{\csp}{{\mathcal{C}\mathcal{P}}}
\newcommand{\VgD}{{\mathcal{V}_{g_D}}}
\newcommand{\ldb}{\left\llbracket}
\newcommand{\rdb}{\right\rrbracket}
\newcommand{\mcv}{{H}}
\newcommand{\ct}{\mcH_\Sigma}
\newcommand{\Sigmah}{{\Sigma_h}}
\newcommand{\nablas}{\nabla_\Sigma}
\newcommand{\nablash}{\nabla_{\Sigma_h}}
\newcommand{\rhos}{\rho_\Sigma}
\newcommand{\Ds}{\text{D}_\Sigma}
\newcommand{\divs}{\text{div}_\Sigma}
\newcommand{\lp}{\left(}
\newcommand{\rp}{\right)}

\newcommand{\tn}{|\mspace{-1mu}|\mspace{-1mu}|}
\newcommand{\curl}{{\text{\rm curl}\,}}

\newcommand{\Vnnsym}{[V_h]^{n\times n}_{sym}}
\global\long\def\bf#1{\mathbf{#1}}

\global\long\def\tr#1{#1^{\intercal}}

\title{Finite element procedures for computing normals and mean curvature
on triangulated surfaces and their use for mesh refinement}

\author[1]{Mirza Cenanovic}

\author[1]{Peter Hansbo}

\author[2]{Mats G. Larson}

\affil[1]{Department of Mechanical Engineering, J\"{o}nk\"{o}ping University, SE-55111
J\"{o}nk\"{o}ping, Sweden}

\affil[2]{Department of Mathematics and Mathematical Statistics, {Ume\aa} University,
SE-901 87 {Ume\aa}, Sweden}

\maketitle
\begin{abstract}
In this paper we consider finite element approaches to computing the
 mean curvature vector and normal at the vertices of piecewise linear triangulated surfaces.
In particular, we adopt a stabilization technique which allows for first order $L^2$--convergence
of the mean curvature vector
and apply this stabilization technique also to the computation of continuous,
recovered, normals using $L^2$--projections of the piecewise
constant face normals. Finally, we use our projected normals
to define an adaptive mesh refinement approach to geometry resolution where we also
employ spline techniques to reconstruct the surface before refinement.
We compare or results to previously proposed approaches.
\end{abstract}



\section{Introduction}

Our aim in this paper is to apply finite element techniques for computing geometrical quantities of interest in 
computer graphics applications, and to show that they can give accurate results, indeed more accurate that classical approaches. 
We restrict ourselves to closed surfaces approximated by piecewise linear simplices, and on such surfaces
we consider three issues:
\begin{itemize}
\item accurate computation of the mean curvature vector;
\item accurate computation of surface normals;
\item adaptive refinement techniques
for resolving the curvature.
\end{itemize}
We discretize the normal and curvature vectors using a piecewise linear finite element method based on tangential differential calculus, following
the approach initiated by Dziuk \cite{Dz88}. This results in piecewise linear, continuous, vector fields on the discrete surface. In order to make comparisons
with standard methods of computing curvature and normals, which are typically only represented at the vertices of the triangulated surface, we will focus mainly on
the nodal values of the finite element fields.
\newline\newline
\textbf{Mean curvature.} The mean curvature vector on a discrete surface
plays an important role in computer graphics and computational geometry,
as well as in certain surface evolution problems, see, e.g. \cite{BoKoPaAlLe10,BoSo08,CeHaLa14,DeMeScBa99,Dz91,Dz08,Sch96}.
It can be obtained by letting the Laplace\textendash Beltrami
operator act on the embedding of the surface in $\mathbb{R}^{3}$,
and various formulas based on this fact have been suggested in the literature, see \cite{MeyDesSchBa03}
and the references therein. It is known that the standard mean curvature
vector based on the finite element discretization of the Laplace\textendash Beltrami
operator on a piecewise linear triangulated surface cannot be expected, in general, to give any order of convergence
in the $L^{2}$ norm. More generally, for triangulated piecewise polynomial
surfaces of order $k$ the expected convergence in $L^{2}$ norm is $k-1$,
cf. \cite{Hei06,De09}. Convergence will also not occur in
other standard discretization methods without
restrictive assumptions on the mesh, see \cite{Xu13}. 
In this paper we employ a stabilized piecewise linear finite element method first suggested in \cite{HaLaZa15} for approximation
of the mean curvature vector, giving first order convergence in the $L^2$ norm for piecewise linear surfaces. The stabilization consists of adding
suitably scaled terms involving the jumps in the tangent gradient
of the discrete mean curvature vector in the direction of the outer
co-normals at each edge in the surface mesh to the $L^{2}$--projection
of the discrete Laplace\textendash Beltrami
operator used to compute the discrete mean curvature vector. 
\newline\newline
\textbf{Normal vectors.} Accurately determining the vertex normals on triangulated surfaces
is of great importance in computer graphics for the computation of
smooth shading \cite{Gouraud1971,Phong1975}, and it is important in surface
meshing/re-meshing \cite{Surazhsky2003,Neto2013,Vlachos2001} as well
as smoothing (fairing) techniques \cite{Hildebrandt2004}. 
We here extend the method suggested for computing the mean curvature vector, which can be seen as a general stabilization approach, 
to the problem of computing accurate vertex normals by stabilized $L^2$--projections. 
\newline\newline
\textbf{Adaptive mesh improvement.} Mesh improvement when the geometry is given by an analytical expression (or is otherwise known)
can be obtained by local refinement of the simplices, putting new vertices on the known surface. The goal is then to resolve the curvature of the mesh in some predefined way.
We suggest an approach based on the difference between the piecewise constant facet normals and the computed finite element normal field.
This gives an estimate of the error in discrete facet normals which is closely related to the curvature of the geometry as will be discussed below.
If the geometry is not a priori known but we are simply given a point cloud or a mesh, interpolation using vertex normals is standard, cf., e.g., Boschiroli
et al. \cite{Boschiroli2011}. We combine one such approach, the PN triangle of Vlachos et al. \cite{Vlachos2001}, with our finite element normal fields and adaptive scheme in order to 
enhance the refined geometry.
\newline

The outline of the remainder of the paper is as follows: In Section
2 we introduce the discrete surface approximations, in Section 3 we
define the stabilized mean curvature vector, in Section 4 we discuss a
different schemes for computing vertex normals, including our stabilized projection method,
in Section 5 we present an adaptive algorithm for resolving curvature, and
in Section 6 we give some representative numerical results.

\section{Meshed surfaces}

Consider an embedded orientable closed surface $\mathbb{R}^{3}\supset\Sigma\in C^{2}$
with exterior unit normal $\bm{n}$. Let $\phi$ be the signed distance
function such that $\nabla\phi=\bm{n}$ on $\Sigma$ and let $\bm{p}(\bm{x})=\bm{x}-\phi(\bm{x})\bm{n}(\bm{p}(\bm{x}))$
be the closest point mapping. Let $U_{\delta}(\Sigma)$ be the open
tubular neighborhood 
\[
U_{\delta}(\Sigma)=\{\bm{x}\in\mathbb{R}^{3}:|\phi(\bm{x})|<\delta\}
\]
for $\delta>0$ of $\Sigma$. Then there is $\delta_{0}>0$ such that
the closest point mapping $\bm{p}(\bm{x})$ assigns precisely one
point on $\Sigma$ to each $\bm{x}\in U_{\delta_{0}}(\Sigma)$. 

We triangulate $\Sigma$ using a elementwise planar mesh $\mathcal{K}_{h}$
to obtain a quasiuniform triangulated surface 
\[
\Sigma_{h}=\cup_{K\in\mathcal{K}_{h}}K\subset U_{\delta_{0}}(\Sigma).
\]

Using the closest point mapping any function $v$ on $\Sigma$
can be extended to $U_{\delta_{0}}(\Sigma)$ using the pull back 
\begin{equation}
v^{e}=v\circ\bm{p}\quad\text{on \ensuremath{U_{\delta_{0}}(\Sigma)}}
\end{equation}
and the lifting $w^{l}$ of a function $w$ defined on $\Sigma_{h}$
to $\Sigma$ is defined as the push forward 
\begin{equation}
(w^{l})^{e}=w^{l}\circ\bm{p}=w\quad\text{on \ensuremath{\Sigma_{h}}}\label{eq:wle}
\end{equation}

\section{Approximation of the mean curvature vector}

\subsection{The continuous mean curvature vector}

We define the tangential surface gradient $\nabla_{\Sigma}$ by $\nabla_{\Sigma}:=\bm{P}_{\Sigma}\nabla$,
where $\nabla$ is the $\mathbb{R}^{3}$ gradient and $\bm{P}_{\Sigma}(\bm{x})=\bm{I}-\bm{n}(\bm{x})\otimes\bm{n}(\bm{x})$
is the projection onto the tangent plane $T_{\Sigma}(\bm{x})$ of
$\Sigma$ at point $\bm{x}\in\Sigma$. The mean curvature vector $\bm{H}:\Sigma\rightarrow\mathbb{R}$
is then defined by 
\begin{equation}
\bm{H}=-\Delta_{\Sigma}\bm{x}_{\Sigma}
\end{equation}
where $\bm{x}_{\Sigma}:\Sigma\ni\bm{x}\mapsto\bm{x}\in\mathbb{R}^{3}$
is the coordinate map of $\Sigma$ into $\mathbb{R}^{3}$ and $\Delta_{\Sigma}=\nabla_{\Sigma}\cdot\nabla_{\Sigma}$
is the Laplace\textendash Beltrami operator. 

The relation between the mean curvature vector and mean curvature
is given by the identity 
\begin{equation}
\bm{H}=(\kappa_{1}+\kappa_{2})\bm{n}
\end{equation}
where $\kappa_{1}$ and $\kappa_{2}$ are the two principal curvatures
and $(\kappa_{1}+\kappa_{2})/2=:H$ is the mean curvature, see \cite{BoKoPaAlLe10}.

The mean curvature vector satisfies the following weak problem: find
$\bm{H}\in W=[H^{1}(\Sigma)]^{3}$ such that 
\begin{equation}
(\bm{H},\bm{v})_{\Sigma}=(\nabla_{\Sigma}\bm{x}_{\Sigma},\nabla_{\Sigma}\bm{v})_{\Sigma}\quad\forall\bm{v}\in W\label{meancurvature}
\end{equation}
where $\nabla_{\Sigma}\bm{w}=\bm{w}\otimes\nabla_{\Sigma}$ for a
vector valued function $\bm{w}$ and 
\[
(\bm{v},\bm{w})_{\omega}=\int_{\omega}\bm{v}\cdot\bm{w}dx
\]
is the $L^{2}$--inner product on the set $\omega$ with associated
norm 
\[
\|\bm{v}\|_{\omega}^{2}=\int_{\omega}\bm{v}\cdot\bm{v}dx.
\]

\begin{figure}
\begin{centering}
\includegraphics[width=0.5\textwidth]{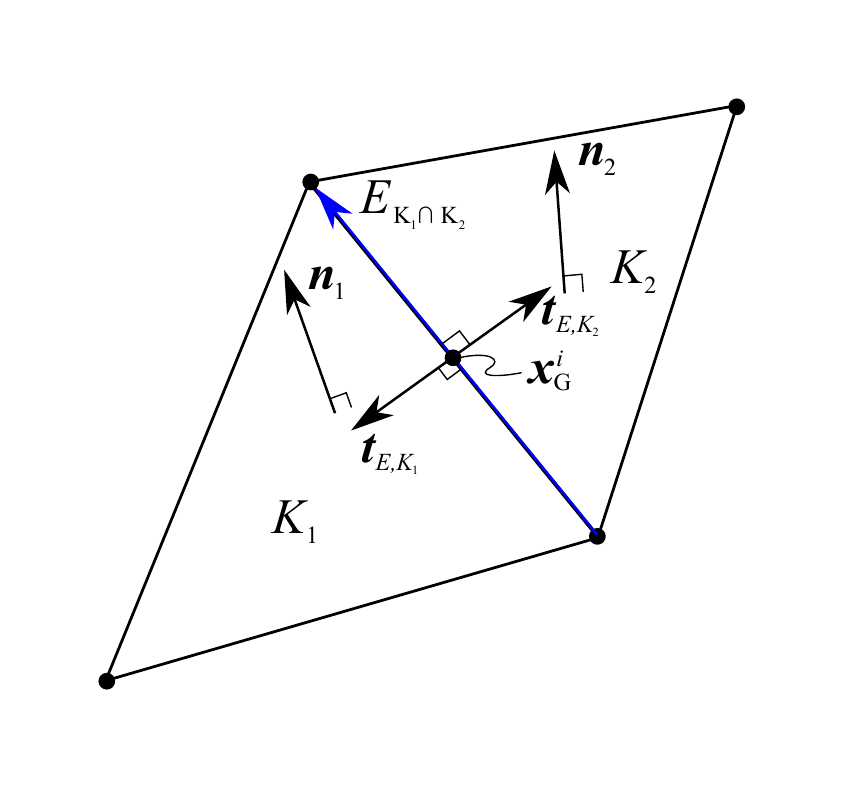}\caption{Edge stabilization\label{fig:Edge-stabilization}}
\end{centering}
\end{figure}

Given the discrete coordinate map $\bm{x}_{\Sigma_{h}}:\Sigma_{h}\ni\bm{x}\mapsto\bm{x}\in\mathbb{R}^{3}$
and a discrete projection operator $\bm{P}_{\Sigma_h}=\bm{I}-\bm{n}_h\otimes\bm{n}_n$, where $\bfn_h$ denotes the piecewise constant facet normals, we define the stabilized
discrete mean curvature vector $\bm{H}_{h}$ as follows. Let $V_{h}$ be the space of piecewise linear
continuous functions defined on $\mathcal{K}_{h}$ and seek $\bm{H}_{h}\in W_{h}=[V_{h}]^{3}$
such that 
\begin{equation}
(\bm{H}_{h},\bm{v})_{\Sigma_{h}}+J_{h}(\bm{H}_{h},\bm{v})=(\nabla_{\Sigma_{h}}\bm{x}_{\Sigma_{h}},\nabla_{\Sigma_{h}}\bm{v})_{\Sigma_{h}}\quad\forall\bm{v}\in W_{h}\label{meancurvaturedisc}
\end{equation}
where $\nabla_{\Sigma_{h}}=\bm{P}_{\Sigma_h}\nabla$ and the stabilization
term $J_{h}(\cdot,\cdot)$ is defined by 
\begin{align}
J_{h}(u,v) & =\gamma\sum_{E\in\mathcal{E}_{h}}h([\bm{t}_{E}\cdot\nabla_{\Sigma_{h}}u],[\bm{t}_{E}\cdot\nabla_{\Sigma_{h}}v])_{E}.
\end{align}
Here $\gamma\geq0$ is a stabilization parameter and $\mathcal{E}_{h}=\{E\}$
is the set of edges in the partition $\mathcal{K}_{h}$ of $\Sigma_{h}$.
The jump of the tangential derivative in the direction of the outer
co-normals at an edge $E\in\mathcal{E}_{h}$ shared by elements $K_{1}$
and $K_{2}$ in $\mathcal{K}_{h}$ is defined by 
\begin{equation}
[\bm{t}_{E}\cdot\nabla_{\Sigma_{h}}u]=\bm{t}_{E,K_{1}}\cdot\nabla_{\Sigma_{h}}u_{1}+\bm{t}_{E,K_{2}}\cdot\nabla_{\Sigma_{h}}u_{2}
\end{equation}
where $u_{i}=u|_{K_{i}}$, $i=1,2,$ and $\bm{t}_{E,K_{i}}$ are the
co-normals, i.e., the unit vectors orthogonal to $E$, tangent and
exterior to $K_{i}$, $i=1,2$, see Figure \ref{fig:Edge-stabilization}.
This stabilization method allows for proving first order convergence
of the curvature vector $\|\bm{H}-\bm{H}_{h}^{l}\|_{\Sigma}\lesssim h,$
see \cite{HaLaZa15}.

\subsection{Implementation issues}

Using the standard Galerkin approximation, 
\[
u\approx\sum_{i}\varphi_{i}U_{i},
\]
where $\varphi_{i}$ is the finite element basis functions and $U_{i}$ the
nodal approximations of $u$ we have that 
\[
\nabla_{\Sigma_{h}}u{\approx}\sum_{i}\nabla_{\Sigma_{h}}\varphi_{i}U_{i},
\]
were we define the tangential gradient of the basis function by

\renewcommand*{\arraystretch}{2}

\[
\nabla_{\Sigma_{h}}\varphi_{i}=:\begin{bmatrix}\dfrac{\partial{\varphi_{i}}}{\partial x_{\Sigma_{h}}}\\
\dfrac{\partial{\varphi_{i}}}{\partial y_{\Sigma_{h}}}\\
\dfrac{\partial{\varphi_{i}}}{\partial z_{\Sigma_{h}}}
\end{bmatrix}=\bm{P}_{{\Sigma_h}}\begin{bmatrix}\dfrac{\partial{\varphi_{i}}}{\partial x}\\
\dfrac{\partial{\varphi_{i}}}{\partial y}\\
\dfrac{\partial{\varphi}_{i}}{\partial z}
\end{bmatrix}.
\]

The tangential derivative of the basis function is given by

\[
(\bm{t}\cdot\nabla_{\Sigma_{h}}){\varphi_{i}}=t_{x}\frac{\partial{\varphi_{i}}}{\partial x_{\Sigma_{h}}}+t_{y}\frac{\partial{\varphi_{i}}}{\partial y_{\Sigma_{h}}}+t_{z}\frac{\partial{\varphi_{i}}}{\partial z_{\Sigma_{h}}}.
\]
For vector\textendash valued unknowns $\bm{u}$ we have $\bm{u}\approx\bm{\Phi}{\bf u}$
where ${\bf u}$ denotes nodal values and 
\begin{equation}
\bm{\Phi}:=\begin{bmatrix}\varphi^{1} & 0 & 0 & \varphi^{2} & 0 & 0 & \cdots\\
0 & \varphi^{1} & 0 & 0 & \varphi^{2} & 0 & \cdots\\
0 & 0 & \varphi^{1} & 0 & 0 & \varphi^{2} & \cdots
\end{bmatrix},\label{phimat}
\end{equation}
and using the notation $\bm{t}_{1}$ and $\bm{t}_{2}$ for the two
co-normals on a given edge $E$, we define 
\[
{\bf B_{E}}:=\left[(\bm{t}_{1}\cdot\nabla_{\Sigma_{h}})\bm{\Phi},(\bm{t}_{2}\cdot\nabla_{\Sigma_{h}})\bm{\Phi}\right],
\]
and the discrete stabilization matrix is given by

\[
{\bf J=\sum_{E\in\mathcal{E}}h\int_{E}\tr{{\bf B_{E}}}{\bf B_{E}dE}}.
\]

The linear system corresponding to \eqref{meancurvaturedisc} becomes

\begin{equation}\label{stabilizedcurvature}
\left({\bf M}+\gamma_{H}{\bf J}\right){\bf H}={\bf S}{\bf x},
\end{equation}
where ${\bf M}$ is the so called mass matrix, given by 
\[
{\bf M=\int_{\Sigma_{h}}\tr{\bm{\Phi}}\bm{\Phi}\,dx,}
\]
$\gamma_{H}$ is the mean curvature specific stabilization factor,
${\bf S}$ is the discrete Laplace-Beltrami operator defined by

\[
{\bf S=\int_{\Sigma_{h}}\tr{\left(\nabla_{\Sigma_{h}}\bm{\varphi}\right)}\nabla_{\Sigma_{h}}\bm{\varphi}\,dx=\int_{\Sigma_{h}}\tr{{\bf B_{S}}}{\bf B_{S}\ dx,}}
\]
with

\renewcommand*{\arraystretch}{2}

\[
{\bf B_{S}:=\begin{bmatrix}\dfrac{\partial\varphi_{1}}{\partial x_{\Sigma_{h}}} & 0 & 0 & \dfrac{\partial\varphi_{2}}{\partial x_{\Sigma_{h}}} & 0 & 0 & \cdots\\
\dfrac{\partial\varphi_{2}}{\partial y_{\Sigma_{h}}} & 0 & 0 & \dfrac{\partial\varphi_{2}}{\partial y_{\Sigma_{h}}} & 0 & 0 & \cdots\\
\dfrac{\partial\varphi_{1}}{\partial z_{\Sigma_{h}}} & 0 & 0 & \dfrac{\partial\varphi_{2}}{\partial z_{\Sigma_{h}}} & 0 & 0 & \cdots\\
0 & \dfrac{\partial\varphi_{1}}{\partial x_{\Sigma_{h}}} & 0 & 0 & \dfrac{\partial\varphi_{2}}{\partial x_{\Sigma_{h}}} & 0 & \cdots\\
0 & \dfrac{\partial\varphi_{1}}{\partial y_{\Sigma_{h}}} & 0 & 0 & \dfrac{\partial\varphi_{2}}{\partial y_{\Sigma_{h}}} & 0 & \cdots\\
\vdots & \vdots & \vdots & \vdots & \vdots & \vdots & \ddots
\end{bmatrix},}
\]
\renewcommand*{\arraystretch}{1}${\bf x}$ is the coordinate vector
of the nodal positions in the mesh, and ${\bf H}$ denotes the vector
of vertex values of the approximate mean curvature vector.

\subsection{Alternative approximations of the mean curvature vector\label{subsec:DiscreteLocalLaplaceBeltrami}}

There exists several well known approaches to mean curvature estimation, for
an extensive overview, see \cite{Magid2007}. In the context of finite elements
an alternative to ours is proposed by Heine in \cite{Heine2004}.

\begin{figure}
\begin{centering}
\subfloat[\label{fig:1-ring-neighborhood}]{\centering{}\includegraphics[width=0.35\textwidth]{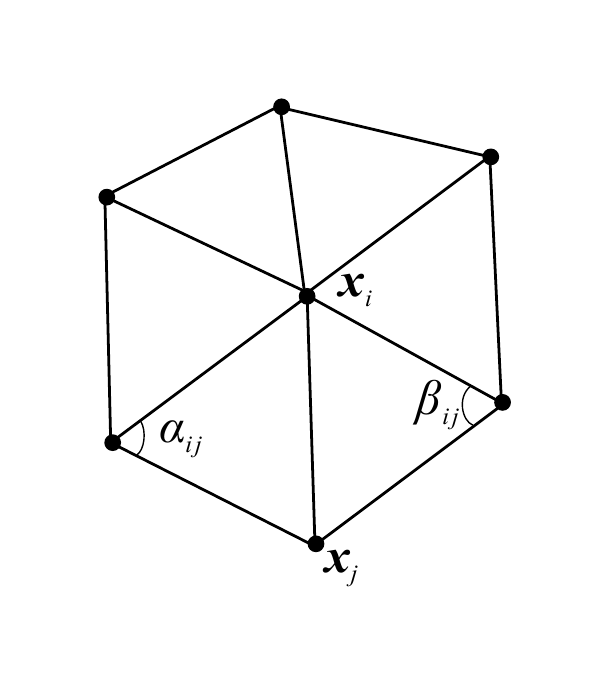}}\subfloat[\label{fig:1-ring-neighborhood-1}]{\centering{}\includegraphics[width=0.35\textwidth]{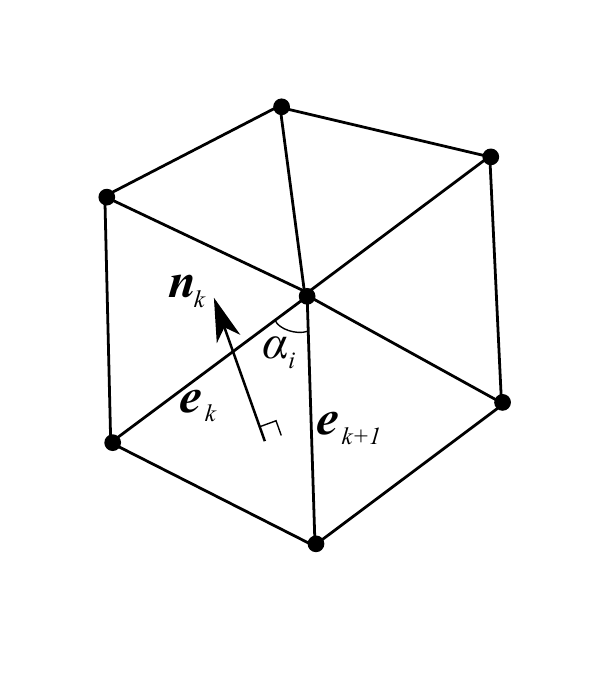}}\caption{1-ring neighborhood to ${\bf x_{i}}$}
\par\end{centering}
\end{figure}

We shall compare
our method to two types of such approaches:
1) fitting a surface locally to each vertex, see e.g., \cite[Chap. 8.5]{Berentzen2012}
and 2) computing the discrete local Laplace-Beltrami operator, see
e.g., \cite{MeyDesSchBa03,DeMeScBa99}. 
\newline\newline
\textbf{Smooth surface fit.}
Curvatures can be computed using a locally fitted quadratic function
\[
f(u,v)=\dfrac{1}{2}\left(au^{2}+2buv+cv^{2}\right)
\]
around a point ${\bf x_{i}}$ with $u$ and $v$ are local coordinates
of the tangential plane to ${\bf x_{i}}$ such that $f(0,0)=0$. The
tangential plane is determined using one of the edges connected to
${\bf x_{i}}$ and the normal at the same point. The idea is to compute
the shape operator or Weingarten map of this function and subsequently
the curvature. See \cite[Chap. 8.5]{Berentzen2012} for further details.
\newline\newline
\textbf{The discrete local Laplace-Beltrami operator.}
Let ${\bf K}$ denote the Laplace\textendash Beltrami operator so
that ${\bf K}(\bm{x})=2H(\bm{x})\bm{n}(\bm{x})$ at a given point
$\bm{x}$ on the surface. On triangulated surfaces, one can use Gauss'
theorem to extract a discrete version of this operator in the nodes
${\bf x_{i}}$ of the mesh, cf. Meyer et al. \cite{MeyDesSchBa03}.
The integral of ${\bf K}$ over the discrete 1-ring surface M on a
triangulated surface is then given by 
\[
\int_{A_{M}}{\bf K({\bm{x})dA=\dfrac{1}{2}\sum_{j\in N}\left(\cot\alpha_{ij}+\cot\beta_{ij}\right)\left({\bf x_{i}-{\bf x_{j}}}\right),}}
\]
where the angles $\alpha_{ij}$ and $\beta_{ij}$ are opposite to
the edge $i\ j$ and $N$ is the set of neighbour vertices to ${\bf x_{i}}$,
see Figure \ref{fig:1-ring-neighborhood}. Given some definition $A_{V}$
of the local area surrounding a vertex ${\bf x_{i}}$ we can then
define the discrete approximation ${\bf K}_h$ of ${\bf K}$ as 
\[
{\bf K}_{h}({\bf x_{i}}):=\dfrac{1}{2A_{V}}\sum_{j\in N}\left(\cot\alpha_{ij}+\cot\beta_{ij}\right)\left({\bf x_{i}-{\bf x_{j}}}\right).
\]
In \cite{MeyDesSchBa03}, it is proposed to use the Voroni regions
as the definition for the local area, and an algorithm to improve
the robustness for arbitrary meshes was provided. Similarly, Desbrun
et. al \cite{DeMeScBa99} used the barycentric area to average the
discrete Laplacian. In both cases, in order to compute the vertex
normal, ${\bf K_{h}({\bf x_{i})}}$ is simply normalized and cases
where the curvature is zero are treated by computing the mean face-normal
of the 1-ring neighbourhood. The mean (discrete) curvature at the
vertices is then given by

\[
H_{h}({\bf x_{i}})=\frac{1}{2}\vert{\bf K_{h}({\bf x_{i})\vert.}}
\]

\section{Normal vector approximation}

\subsection{Stabilized projection of the normal field}

In analogy with (\ref{meancurvaturedisc}) we define the recovered
discrete normal vector $\bm{n}_{h}$ as follows: find $\bm{n}_{h}\in W_{h}=[V_{h}]^{3}$
such that 
\begin{equation}
(\bm{n}_{h},\bm{v})_{\Sigma_{h}}+J_{h}(\bm{n}_{h},\bm{v})=(\bm{n}_{K},\bm{v})_{\Sigma_{h}}\quad\forall\bm{v}\in W_{h}\label{eq:stabVertexNormal}
\end{equation}
where $\bm{n}_{K}$ is the piecewise constant exterior normal to the
facet elements $K$. The corresponding linear system becomes

\begin{equation}
\left({\bf M+\gamma_{n}{\bf J}}\right){\bf n}_{h}={\bf b,}\label{eq:discreteStabVertexNormal}
\end{equation}
where $\gamma_{n}$ is the normal-specific stabilization factor, ${\bf n}_{h}$ the vector of vertex normals, and
\[
{\bf b}:=\int_{\Sigma_{h}}\tr{\bm{\Phi}}\bm{n}_{K}\,d\Sigma_h.
\]
Note that (\ref{eq:discreteStabVertexNormal}) can be efficiently
solved using a conjugate gradient method since ${\bf M}$ is symmetric,
positive definite and sparse. 

When translating the computed normal vector field to a set of discrete vertex normals, these will here be normalized (the nodal vectors contained in ${\bf n}_{h}$ are not
in general of unit length). 

\subsection{Alternative approaches to computing vertex normals}

Traditionally, vertex normals are estimated either from a local neighborhood
of surrounding face normals using some type of local averaging, see
e.g., \cite{Jin2005,Neto2013} and the references therein. Other estimation
methodologies also exists such as local smooth surface fits, see,
e. g., \cite{Meek2000}. We use the notations for the local vertex
normals introduced in \cite{Neto2013} and give a brief description;
see Figure \ref{fig:1-ring-neighborhood-1} for an explanation of
the notations used.
\newline\newline
\textbf{Mean weighted equally.} Arguably, the most widespread estimation of the vertex normal was
introduced by Gouraud \cite{Gouraud1971} as 
\begin{equation}
\bm{n}_{\text{MWE}}:=\frac{\sum_{i=1}^{n}\bm{n}_{i}}{\left\vert \sum_{i=1}^{n}\bm{n}_{i}\right\vert },
\end{equation}
where $\bm{n}_{i}$ is the face-normal of triangle $i$ , $n$ is
the total number of triangles that share a common vertex for which
the vertex normal is to be estimated and $|.|$ denotes the norm.
Note that we shall subsequently omit making the normalization step
of the vertex normal explicit and assume $\bm{n}=\hat{\bm{n}}:=\bm{n}/|\bm{n}|$.
\newline\newline
\textbf{Mean weighted by angle.}
A vertex normal approximation using angles between the inner edges
was proposed by Th\"urrner and W\"uthrich \cite{Thuerrner1998}.

\begin{equation}
\bm{n}_{\text{MWA}}:=\sum_{i=1}^{n}\bm{n}_{i}\alpha_{i},
\end{equation}
where $\alpha_{i}$ is the angle between two edges $\bm{e}_{k}$ and
$\bm{e}_{k+1}$ of a face $i$ sharing the vertex. 
\newline\newline
\textbf{Mean weighted by sine and edge length reciprocals.}
Max \cite{Max1999} proposed several methods of weighting the face
normals, one of which is to weight by the sine and edge length reciprocals
to take into account the difference in lengths of surrounding edges.

\begin{equation}
\bm{n}_{\text{MWSELR}}:=\sum_{i=1}^{n}\frac{\bm{n}_{i}\ \sin\alpha_{i}}{|\bm{e}_{k}|\ |\bm{e}_{k+1}|}.\label{eq:MWSELR}
\end{equation}
\newline\newline
\textbf{Mean weighted by areas of adjacent triangles.}
Another method proposed by Max \cite{Max1999} is to weight the normals
by the area of the face.

\begin{equation}
\bm{n}_{\text{MWAAT}}:=\sum_{i=1}^{n}\bm{n}_{i}|\bm{e}_{k}|\ |\bm{e}_{k+1}|\ \sin\alpha_{i}=\sum_{i=1}^{n}\bm{n}_{i}|\bm{e}_{k}\times\bm{e}_{k+1}|,
\end{equation}
where the symbol $\times$ denotes the vector cross product.
\newline\newline
\textbf{Mean weighted by edge length reciprocals.}
Max \cite{Max1999} also proposed to just use the edge length reciprocals
as weights.

\begin{equation}
\bm{n}_{\text{MWELR}}:=\sum_{i=1}^{n}\frac{\bm{n}_{i}}{|\bm{e}_{k}|\ |\bm{e}_{k+1}|}.
\end{equation}
\newline\newline
\textbf{Mean weighted by square root of edge length reciprocals.}
Finally, Max \cite{Max1999} also suggested to use the square root
of the length reciprocals.

\begin{equation}
\bm{n}_{\text{MWRELR}}:=\sum_{i=1}^{n}\frac{\bm{n}_{i}}{\sqrt{|\bm{e}_{k}|\ |\bm{e}_{k+1}|}}.
\end{equation}
\newline\newline
\textbf{Normal from the discretized local Laplace-Beltrami operator.}
Another approach is to define the normal using the discretized local
Laplace-Beltrami operator (DLLB) defined in Section \ref{subsec:DiscreteLocalLaplaceBeltrami}.
The normal is defined by normalizing the discrete mean curvature vector.
\begin{equation}
\bm{n}_{\text{DLLB}}:={\bf K}/|{\bf K}|.
\end{equation}

In the numerical example below, Section \ref{sec:numex}, we compare
the accuracy of these different approaches.

\section{Adaptive algorithm}

\subsection{Error estimate}

We base our adaptive algorithm on the Zienkiewicz\textendash Zhu approach
\cite{ZiZh87} which employs the difference between recovered derivatives
and actual discrete piecewise derivatives of a finite element solution. By analogy we consider the piecewise constant
normals to play the role of the piecewise derivatives, and compare
these to the $L^2-$projected normals. 

Since we are focusing on vertex normals, and since we will in the following compare methods that only produce such normals, we define a norm which is an approximation of the $L^2$--norm,
\begin{equation}\label{AppL2}
\| \bfn \|_{L_h^2} := \left( \sum_{K\in\mathcal{K}_{h}}\frac13\text{meas}(K)\sum_{i=1}^3 \vert \bfn (\bfx^i_K)\vert^2\right)^{1/2}
\end{equation}
where $\text{meas}(K)$ denotes the area of $K$ and $\bfx^i_K$ the vertex coordinates on $K$. This represents a Newton--Cotes numerical integration scheme for the $L^2(\Sigma_h)$--norm using the vertices as
integration points. 
The error in normals is thus approximated
\[
\|\bm{n}^{e}-\bm{n}_{K}\|_{L^2_{h}}\approx\|\bm{n}_{h}-\bm{n}_{K}\|_{L^2_{h}}
\]
and we aim at achieving
\[
\|\bm{n}_{h}-\bm{n}_{K}\|_{L^2_{h}}\leq\text{TOL}
\]
where TOL is a given tolerance. We note that we also have
\[
\|\bm{n}^{e}-\bm{n}_{K}\|_{L^2_{h}}\approx\|h\nabla_{\Sigma}\bm{n}\|_{\Sigma}
\]
where $h$ is the local mesh size and $\nabla_{\Sigma}\bm{n}$ is the curvature tensor, 
which indicates that we counter large curvature by reduced mesh size for resolution of the geometry.

\subsection{Triangle refinement}

In cases where the exact geometry is not accessible, we consider triangle 
refinement approaches that utilise vertex normals
for interpolation. An overview of such methods is given by Boschiroli
et al. in \cite{Boschiroli2011}. Nagata \cite{Nagata2005} proposed
a simple quadratic interpolation of triangles using vertex normals
and positions at the end-nodes. The approach by Nagata depends on
a curvature parameter that fixes a curvature coefficient in order
to stabilize the method. The curvature coefficient is highly dependent
on the vertex normal, and in cases where normals are near parallel,
the method cannot capture inflections and without a stabilizing parameter,
cusps will be introduced to the surface, see \cite{Neto2013} where
the authors point out this problem and suggest a possible solution.
The solution suggested in \cite{Neto2013} eliminates the problem
of cusps in the interpolated surface but also eliminates the inflection,
since the segment becomes linear. Another approach is to use higher
order interpolation which are able to capture inflection points. 

\subsubsection{PN triangles}

\begin{figure}
\begin{centering}
\includegraphics[width=0.3\columnwidth]{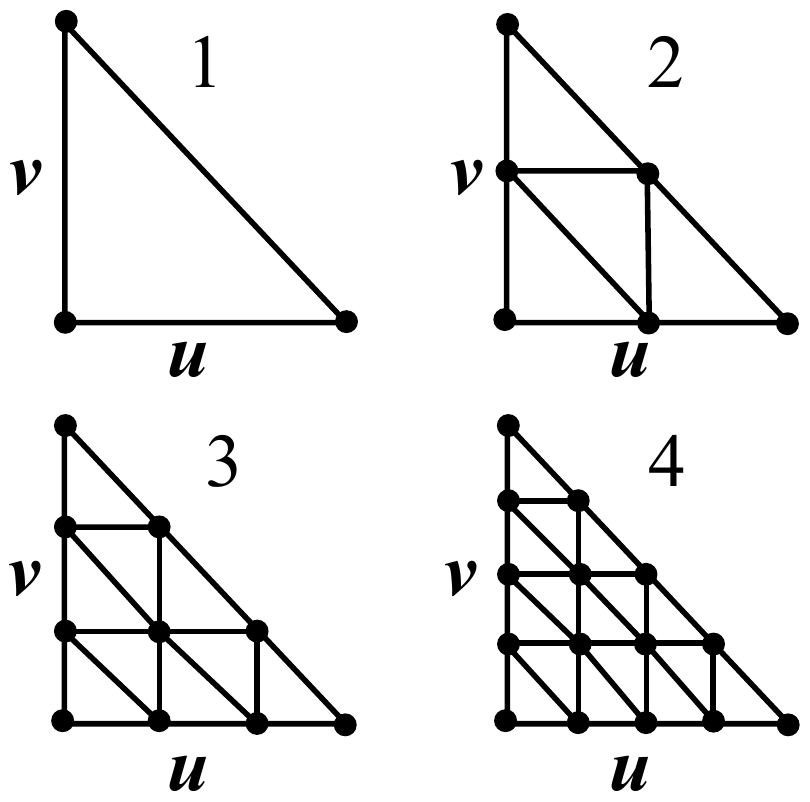}\caption{Subtriangulation created by a uniform $u,v$ grid using $N=\{1,2,3,4\}$. \label{fig:Subtriangulation}}
\par\end{centering}
\end{figure}

\begin{figure}
\begin{centering}
\subfloat[]{
\begin{centering}
\includegraphics[width=0.3\textwidth]{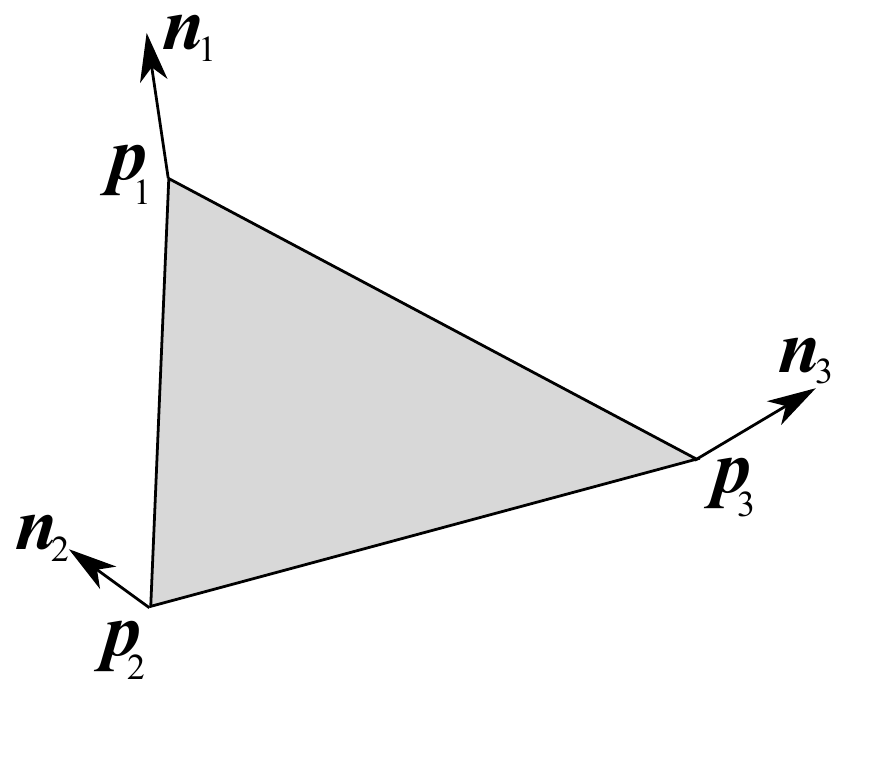}
\par\end{centering}
}\subfloat[]{\centering{}\includegraphics[width=0.3\textwidth]{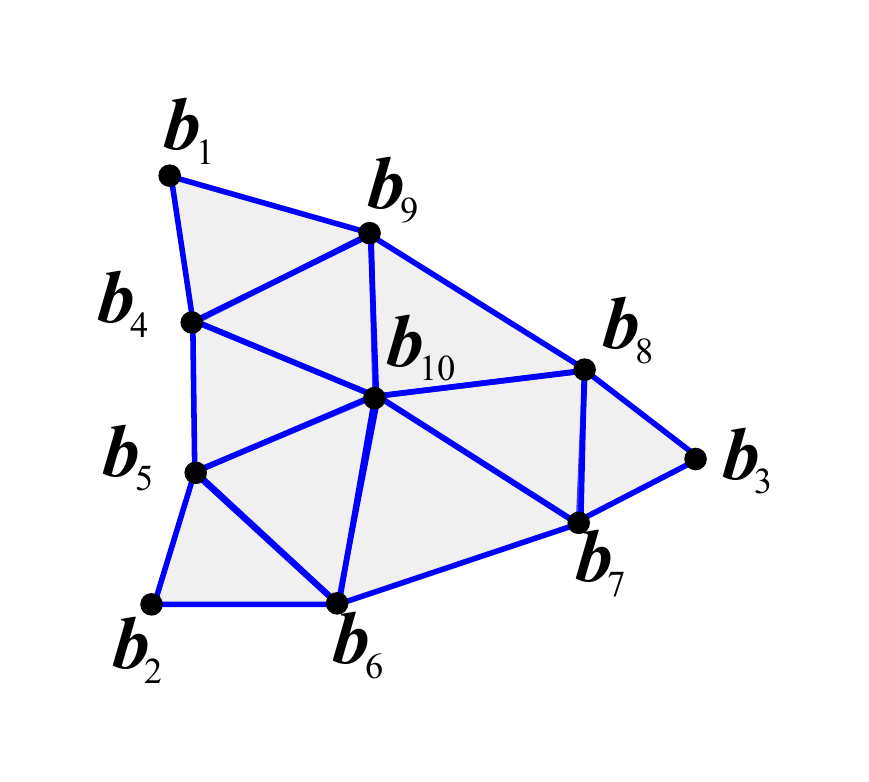}}\caption{PN triangle. a) Initial triangle points $\bm{p}_{i}$ and normals
$\bm{n}_{i}$. b) Control grid with control points $\bm{b}_{k}$.\label{fig:PN-triangle}}
\par\end{centering}
\end{figure}

\begin{figure}
\centering{}\includegraphics[width=0.3\textwidth]{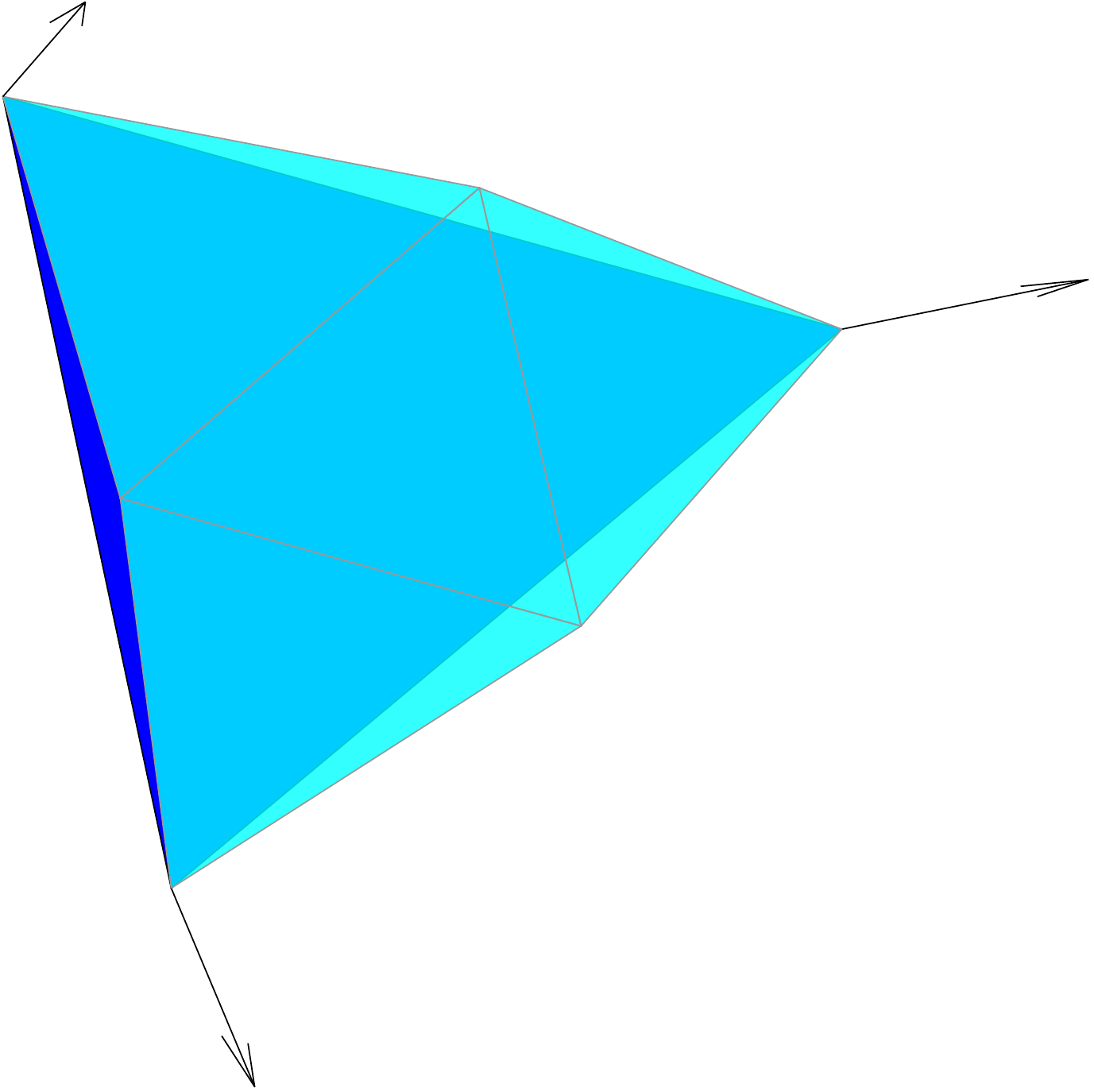}\caption{PN refine algorithm on a triangle using $N=1$ tessellation steps.\label{fig:PNTriangleN1}}
\end{figure}

Vlachos et al. \cite{Vlachos2001} proposed a cubic interpolation
scheme that similarly to Nagata only depends on the positions and
vertex normals of a triangular patch. We here write their algorithm in a vectorized manner. Let then $b:\,\mathbb{R}^{2}\rightarrow\mathbb{\mathbb{R}}^{3}$ denote a cubic triangular
patch given by
\begin{equation}
b(u,v)={\bf B^{\intercal}{\bf U(u,v).}}
\end{equation}
Here ${\bf U}$ is the matrix representation of the parameters defined
by
\begin{equation}
{\bf U=\begin{bmatrix}v^{3} & w^{3} & u^{3} & 3wv^{2} & 3vw^{2} & 3uw^{2} & 3wu^{2} & 3vu & 3uv & 6uvw\end{bmatrix}^{\intercal}}
\end{equation}
where $u=i/N$, $v =j/N$ for $i,j=\{0,1,\ldots,N\}$ such that $w:=1-u-v \geq 0$.
Here $N$ gives a subtriangulation of the initial patch, see Figure \ref{fig:Subtriangulation}.
${\bf B}$ denotes the cubic coefficients in matrix form and is
given by

\begin{equation}
{\bf B=\begin{bmatrix}\bm{b}_{1} & \bm{b}_{2} & \bm{b}_{3} & \bm{b}_{4} & \bm{b}_{5} & \bm{b}_{6} & \bm{b}_{7} & \bm{b}_{8} & \bm{b}_{9} & \bm{b}_{10}\end{bmatrix}^{\intercal}}
\end{equation}
where $\bm{b}$ denote the control points of the control grid for
the PN triangle, see Figure (\ref{fig:PN-triangle}), and are defined
as follows:

\begin{alignat}{1}
\bm{b}_{1} & =\bm{p}_{1},\\
\bm{b}_{2} & =\bm{p}_{2},\\
\bm{b}_{3} & =\bm{p}_{3},\\
\text{for} & \begin{cases}
i=[2,1,3,2,1,3]\\
j=[1,2,2,3,3,1]\\
k=[4,5,6,7,8,9]
\end{cases},\\
w_{ij} & =\left(\bm{p}_{j}-\bm{p}_{i}\right)\cdot\bm{n}_{i},\\
\bm{b}_{k} & =\left(2\bm{p}_{j}+\bm{p}_{i}-w_{ij}\bm{n}_{j}\right)/3,\\
E & =\sum_{k}\bm{b}_{k}/6,\\
V & =\sum_{m=1}^{3}\bm{p}_{m}/3\\
\bm{b}_{10} & =E+\left(E-V\right)/2
\end{alignat}
where $\bm{p}_{i}$ and $\bm{n}_{i}$ are the input corner points
and normals. Finally the total set of interpolated points is given
as a matrix product by

\begin{equation}
{\bf P=\left({\bf B^{\intercal}{\bf U}}\right)^{\intercal}.}
\end{equation}
Note that ${\bf U}$ can be evaluated for a certain number of refinements
$N$ in a pre-processing step. In the local refinement section of
this paper we use $N=1$ see Figure \ref{fig:PNTriangleN1}. As for
the internal vertex normal computation, we do not interpolate the
normals locally, instead we compute $\bm{n}_{h}$ using (\ref{eq:discreteStabVertexNormal})
for the total mesh in each iteration. The reason behind why we limit
the tessellation step to 1 is the subsequent complexity of the local
refinement procedure.

\subsection{Local refinement procedure}

Since the PN refinement with $N=1$ splits the face of a flat triangle
into four child elements, we need a way of handling the hanging nodes.
In this work we adapt the Red-Green refinement method proposed by
Banks et al in \cite{Bank1983}. This method preserves the aspect
ratio of the initial mesh which is crucial in order to secure the
accuracy of the associated finite element method.

\section{Numerical examples\label{sec:numex}}

\subsection{Geometry}

\begin{figure}
\begin{centering}
\subfloat[]{\begin{centering}
\includegraphics[width=0.49\textwidth]{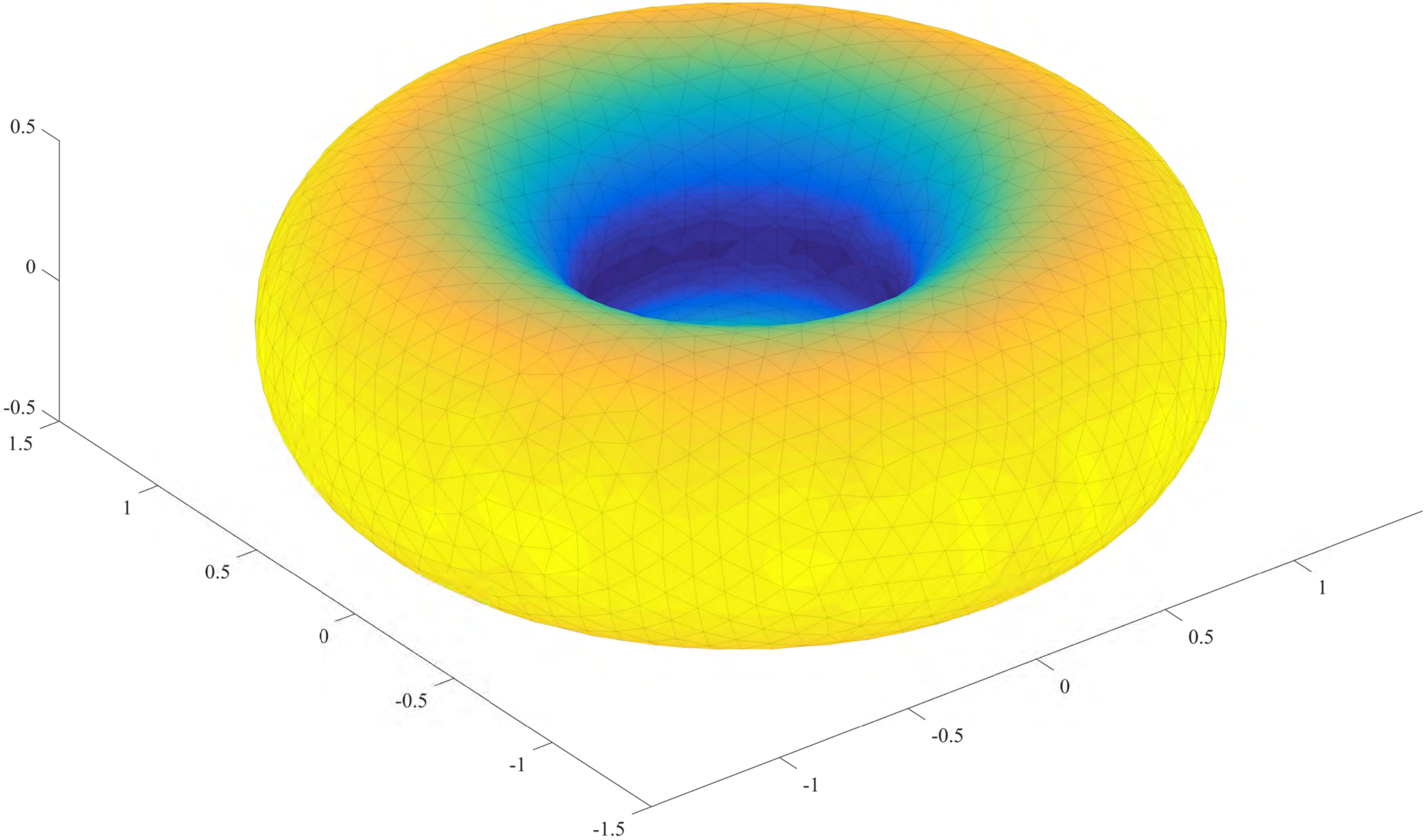}
\par\end{centering}
}$\ $\subfloat[]{\centering{}\includegraphics[width=0.49\textwidth]{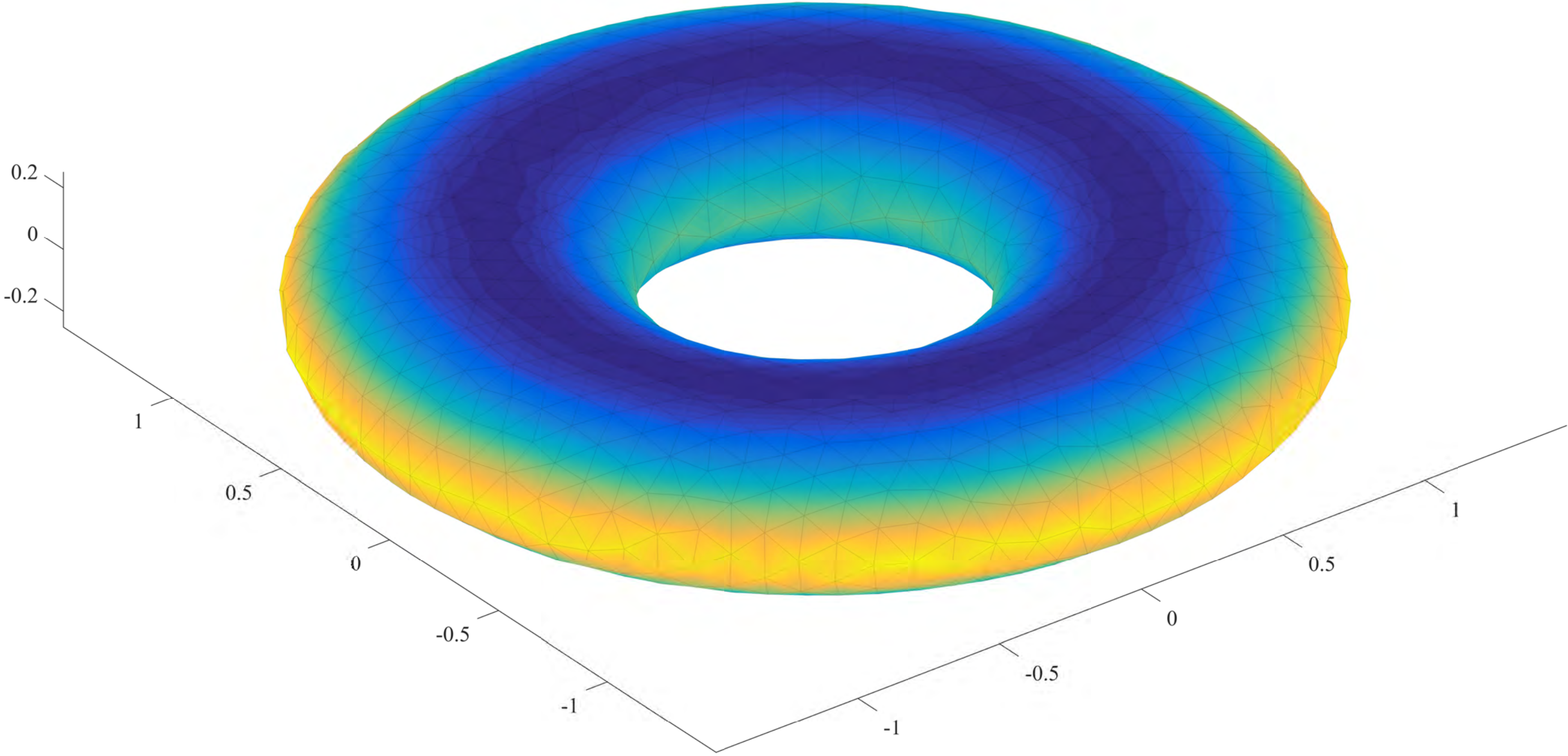}}\caption{Mean curvature of a torus with: (a) $a=1$ and (b) \textbf{$a=4$\label{fig:TorusTypes}}}
\par\end{centering}
\end{figure}

We choose to analyze the errors on an implicitly defined torus which
we can modify in order to generate slightly more complex features.
The surface equation for the torus is given by

\[
\phi(x,y,z)=\left(R-\sqrt{x^{2}+y^{2}}\right)^{2}+az^{2}-r^{2},
\]
where $R$ is the torus radius, $r$ the tube radius and $a$ is a
``squish-factor'' used to squish the torus in the z-direction in
order to induce a higher curvature on the inside and outside, see
Figure \ref{fig:TorusTypes}. In the following, the torus will be analyzed with
$a=1$ and $a=4$, in order to compare errors with respect to strongly
and smoothly varying curvature.

\subsection{Vertex normal error}

What follows is a comparison of different vertex normals with the exact normal.
The measure for the mesh-size used in this context is defined as
\[
h:=\frac{1}{\sqrt{N_\text{v}}},
\]
where $N_\text{v}$ denotes the number of vertices in the mesh. 

Using an implicitly defined surface $\Sigma=\{\bm{x}:\phi(\bm{x})=0\}$,
where $\phi$ is a signed distance function with the property $|\nabla\phi|=1,$
we have that $\bm{n}(\bm{x}_{\Sigma})=\nabla\phi(\bm{x}_{\Sigma})$.
As discussed above, we will use (\ref{AppL2}) and define the error as
\begin{equation}
\epsilon=\|\bm{n}_{a}-\bm{n}_{e}\|_{L^2_h},\label{eq:VertexNormalErrorL2}
\end{equation}
where $\bm{n}_{a}$ is the approximate and $\bm{n}_{e}$ the
exact normal defined by $\bm{n}_{e}=\nabla\phi$, computed at the vertex $i$ using $\bm{n}_{e}^{i}=\nabla\phi(\bm{x}^{i})$. 
The convergence
rates are defined as

\[
p_n=\dfrac{\log(\epsilon_{n+1})-\log(\epsilon_{n})}{\log(h_{n+1})-\log(h_{n})}
\]

\subsection{Evaluation of the accuracy of computed vertex normals}

The vertex normal error analysis was done on an unstructured mesh
of a torus with $R=1$, and $r=1/2$ and $a=\{1,4\}$, see Figure
\ref{fig:TorusTypes}. 

\begin{figure}
\begin{centering}
\subfloat[]{\begin{centering}
\includegraphics[width=0.49\textwidth]{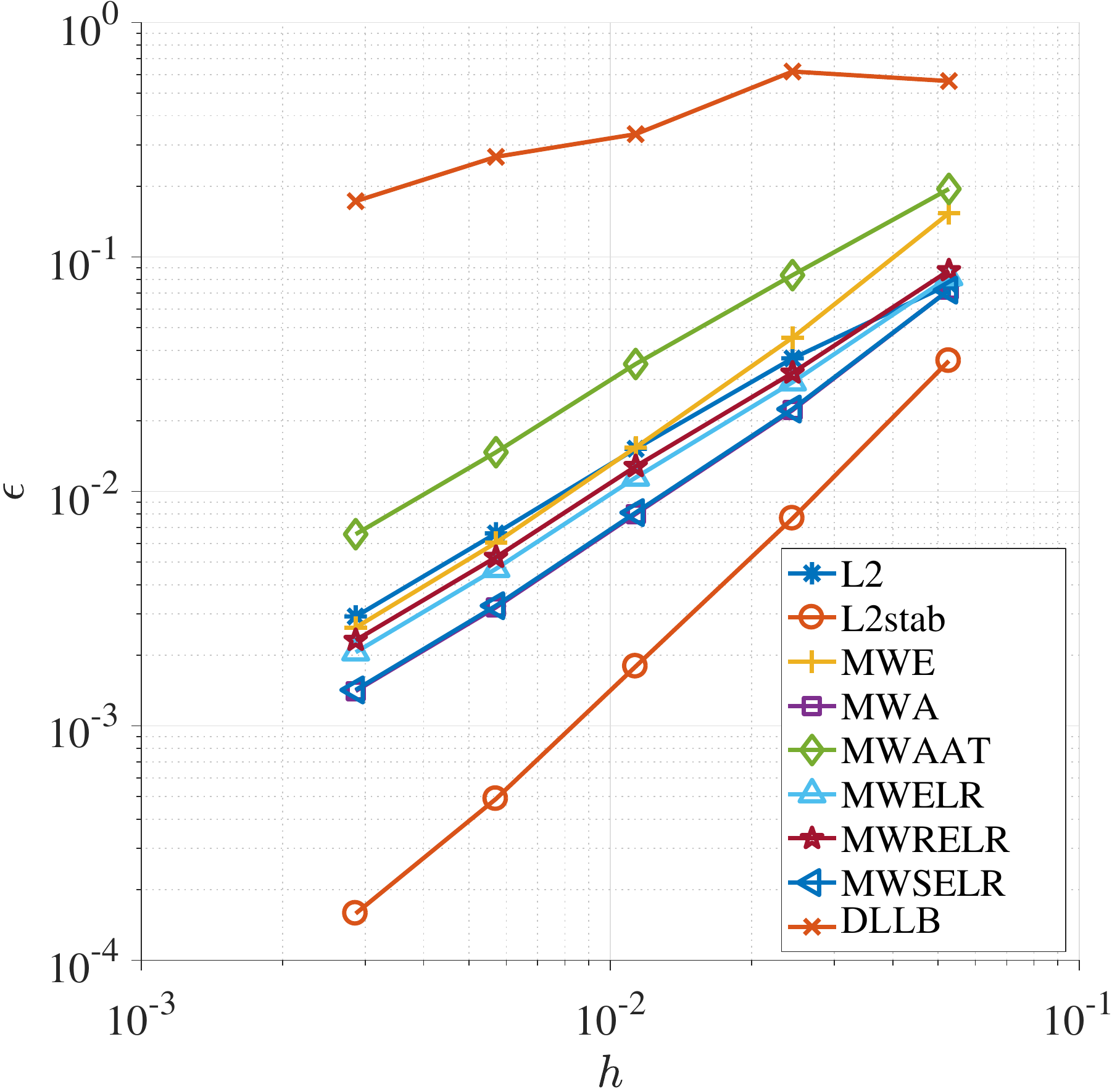}
\par\end{centering}
}$\ $\subfloat[]{\centering{}\includegraphics[width=0.49\textwidth]{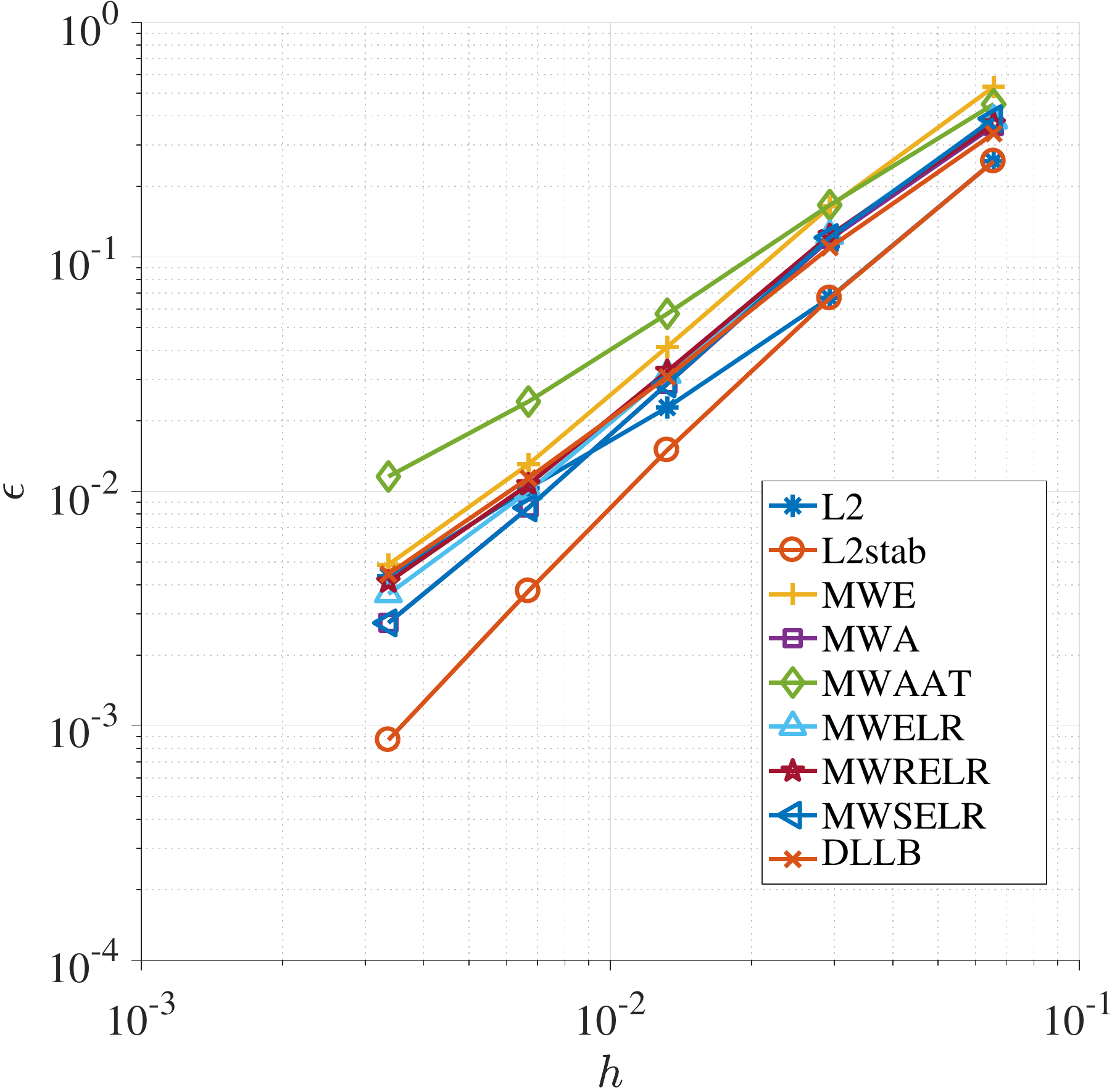}}
\par\end{centering}
\centering{}\caption{Convergence of vertex normals on a torus with: (a) $a=1$ and (b)
$a=4$.\label{fig:ConvNormals}}
\end{figure}

The convergence of $L^2$ errors $\epsilon$ defined in (\ref{eq:VertexNormalErrorL2})
are shown in Figure \ref{fig:ConvNormals} were it can be seen that
the stabilized $L^2$--projection of the normals converges optimally.
The raw data for this graph is available in Table \ref{tab:NormalData}.
The relative difference between the stabilized $L^2$
normals and the next best traditional method $\bm{n}_{\mathrm{MWA}}$
can be seen in Table \ref{tab:DeltaNormalErrors} where we can
see a relative error decrease from $\epsilon_{\text{MWA}}$ of $\sim29\%$
to $\sim88\%$ depending on mesh-size and geometry. The convergence
rates can be viewed in Table \ref{tab:NormaErrlConvRates}. In the
next section we shall analyze the impact of the stabilization on the
normal errors.

\begin{table}
\begin{centering}
\subfloat[Torus with $a=1$]{\begin{centering}
 \resizebox{\textwidth}{!}{
\begin{tabular}{|c|c|c|c|c|c|c|c|c|c|}
\hline 
$h$ & $L^{2}$ & $L_{\mathrm{stab}}^{2}$ & MWE & MWA & MWAAT & MWELR & MWRELR & MWSERL & DLLB\tabularnewline
\hline 
\hline 
0.0527 & 0.0762 & 0.0360 & 0.1535 & 0.0723 & 0.1950 & 0.0820 & 0.0874 & 0.0722 & 0.5623\tabularnewline
\hline 
0.0245 & 0.0369 & 0.0077 & 0.0452 & 0.0222 & 0.0836 & 0.0292 & 0.0320 & 0.0224 & 0.6175\tabularnewline
\hline 
0.0113 & 0.0152 & 0.0018 & 0.0154 & 0.0080 & 0.0349 & 0.0115 & 0.0128 & 0.0081 & 0.3338\tabularnewline
\hline 
0.0057 & 0.0066 & 0.0005 & 0.0060 & 0.0032 & 0.0147 & 0.0047 & 0.0052 & 0.0033 & 0.2672\tabularnewline
\hline 
0.0029 & 0.0029 & 0.0002 & 0.0026 & 0.0014 & 0.0066 & 0.0021 & 0.0023 & 0.0014 & 0.1728\tabularnewline
\hline 
\end{tabular}}
\par\end{centering}}
\par\end{centering}
\subfloat[Torus with $a=4$]{\resizebox{\textwidth}{!}{
\begin{tabular}{|c|c|c|c|c|c|c|c|c|c|}
\hline 
$h$ & $L^{2}$ & $L_{\mathrm{stab}}^{2}$ & MWE & MWA & MWAAT & MWELR & MWRELR & MWSERL & DLLB\tabularnewline
\hline 
\hline 
0.0657 & 0.2557 & 0.2557 & 0.5310 & 0.3618 & 0.4477 & 0.3828 & 0.3701 & 0.3884 & 0.3359\tabularnewline
\hline 
0.0294 & 0.0670 & 0.0666 & 0.1638 & 0.1190 & 0.1664 & 0.1231 & 0.1234 & 0.1209 & 0.1097\tabularnewline
\hline 
0.0132 & 0.0228 & 0.0150 & 0.0413 & 0.0287 & 0.0572 & 0.0314 & 0.0325 & 0.0288 & 0.0308\tabularnewline
\hline 
0.0067 & 0.0103 & 0.0038 & 0.0130 & 0.0085 & 0.0242 & 0.0100 & 0.0107 & 0.0085 & 0.0114\tabularnewline
\hline 
0.0034 & 0.0044 & 0.0009 & 0.0049 & 0.0028 & 0.0116 & 0.0036 & 0.0041 & 0.0027 & 0.0044\tabularnewline
\hline 
\end{tabular}}}\caption{Vertex normal error $\epsilon$ as define in (\ref{eq:VertexNormalErrorL2}).\label{tab:NormalData}}
\end{table}

\begin{table}
\begin{centering}
\subfloat[Torus with $a=1$]{\begin{centering}
\begin{tabular}{|c|c|c|c|c|}
\hline 
$h$ & $\epsilon_{\mathrm{MWA}}-\epsilon_{L^{2}}$ & relative change & $\epsilon_{\mathrm{MWA}}-\epsilon_{L_{\mathrm{stab}}^{2}}$ & relative change\tabularnewline
\hline 
\hline 
0.0527 & -0.0039 & -0.0542 & 0.0363 & 0.5022\tabularnewline
\hline 
0.0245 & -0.0147 & -0.6625 & 0.0145 & 0.6547\tabularnewline
\hline 
0.0113 & -0.0072 & -0.8959 & 0.0062 & 0.7764\tabularnewline
\hline 
0.0057 & -0.0034 & -1.0657 & 0.0027 & 0.8477\tabularnewline
\hline 
0.0029 & -0.0015 & -1.0784 & 0.0012 & 0.8876\tabularnewline
\hline 
\end{tabular}
\par\end{centering}
}
\par\end{centering}
\centering{}\subfloat[Torus with $a=4$]{\begin{centering}
\begin{tabular}{|c|c|c|c|c|}
\hline 
$h$ & $\epsilon_{\mathrm{MWA}}-\epsilon_{L^{2}}$ & relative change & $\epsilon_{\mathrm{MWA}}-\epsilon_{L_{\mathrm{stab}}^{2}}$ & relative change\tabularnewline
\hline 
\hline 
0.0657 & 0.1060 & 0.2931 & 0.1060 & 0.2931\tabularnewline
\hline 
0.0294 & 0.0521 & 0.4375 & 0.0524 & 0.4402\tabularnewline
\hline 
0.0132 & 0.0059 & 0.2064 & 0.0137 & 0.4781\tabularnewline
\hline 
0.0067 & -0.0018 & -0.2078 & 0.0048 & 0.5597\tabularnewline
\hline 
0.0034 & -0.0016 & -0.5827 & 0.0019 & 0.6837\tabularnewline
\hline 
\end{tabular}
\par\end{centering}
}\caption{Error differences, absolute and relative. \label{tab:DeltaNormalErrors} }
\end{table}

\begin{table}
\begin{centering}
\subfloat[Torus with $a=1$.]{\centering{} \resizebox{\textwidth}{!}{
\begin{tabular}{|c|c|c|c|c|c|c|c|c|c|}
\hline 
$h$ & $L^{2}$ & $L_{\mathrm{stab}}^{2}$ & MWE & MWA & MWAAT & MWELR & MWRELR & MWSERL & DLLB\tabularnewline
\hline 
\hline 
0.0527 & - & - & - & - & - & - & - & - & -\tabularnewline
\hline 
0.0245 & 0.9452 & 2.0151 & 1.5909 & 1.5386 & 1.1032 & 1.3449 & 1.3087 & 1.5214 & -0.1219\tabularnewline
\hline 
0.0113 & 1.1524 & 1.8875 & 1.3999 & 1.3230 & 1.1332 & 1.2153 & 1.1924 & 1.3205 & 0.7987\tabularnewline
\hline 
0.0057 & 1.2124 & 1.8968 & 1.3641 & 1.3375 & 1.2636 & 1.3078 & 1.3005 & 1.3333 & 0.3249\tabularnewline
\hline 
0.0029 & 1.1844 & 1.6348 & 1.2098 & 1.1932 & 1.1687 & 1.1897 & 1.1841 & 1.2009 & 0.6328\tabularnewline
\hline 
\end{tabular}}}
\par\end{centering}
\centering{}\subfloat[Torus with $a=4$.]{\centering{} \resizebox{\textwidth}{!}{
\begin{tabular}{|c|c|c|c|c|c|c|c|c|c|}
\hline 
$h$ & $L^{2}$ & $L_{\mathrm{stab}}^{2}$ & MWE & MWA & MWAAT & MWELR & MWRELR & MWSERL & DLLB\tabularnewline
\hline 
\hline 
0.0657 & - & - & - & - & - & - & - & - & -\tabularnewline
\hline 
0.0294 & 1.6651 & 1.6711 & 1.4612 & 1.3812 & 1.2299 & 1.4099 & 1.3647 & 1.4506 & 1.3909\tabularnewline
\hline 
0.0132 & 1.3512 & 1.8705 & 1.7267 & 1.7826 & 1.3383 & 1.7122 & 1.6712 & 1.7999 & 1.5934\tabularnewline
\hline 
0.0067 & 1.1654 & 2.0310 & 1.6914 & 1.7818 & 1.2655 & 1.6789 & 1.6342 & 1.7863 & 1.4589\tabularnewline
\hline 
0.0034 & 1.2520 & 2.1265 & 1.4343 & 1.6452 & 1.0709 & 1.4678 & 1.3909 & 1.6451 & 1.3711\tabularnewline
\hline 
\end{tabular}}}\caption{Vertex normal convergence rates.\label{tab:NormaErrlConvRates}}
\end{table}

\subsection{Effect of the stabilization on the accuracy of the computed normal}

\begin{figure}
\centering{}\includegraphics[width=1\columnwidth]{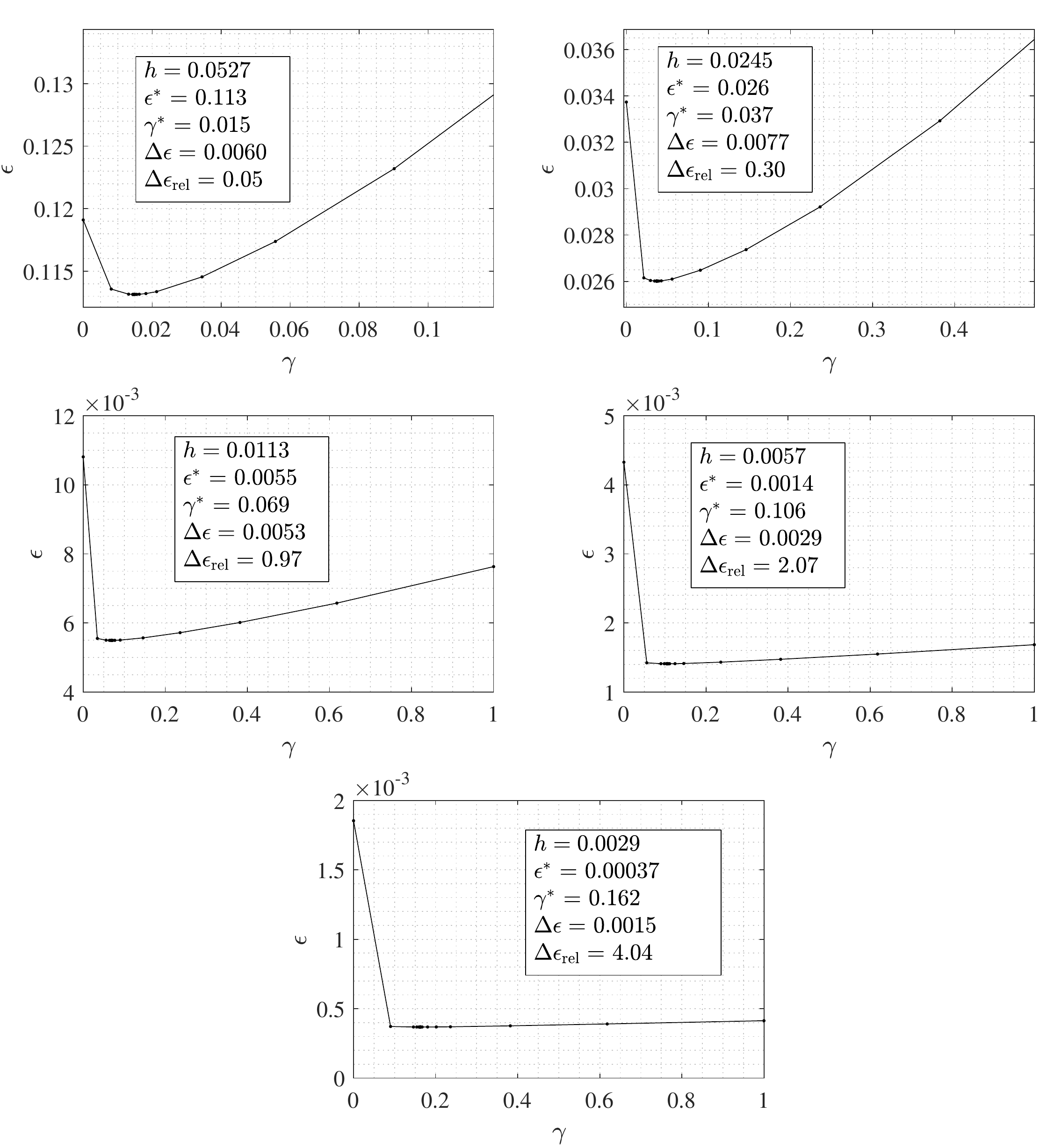}\caption{Vertex normal error as a function of the stabilization factor, $\epsilon(\gamma_{n}),$
for different mesh-sizes on a torus with $a=1$.\label{fig:gammaNormalA1}}
\end{figure}

\begin{figure}
\centering{}\includegraphics[width=1\columnwidth]{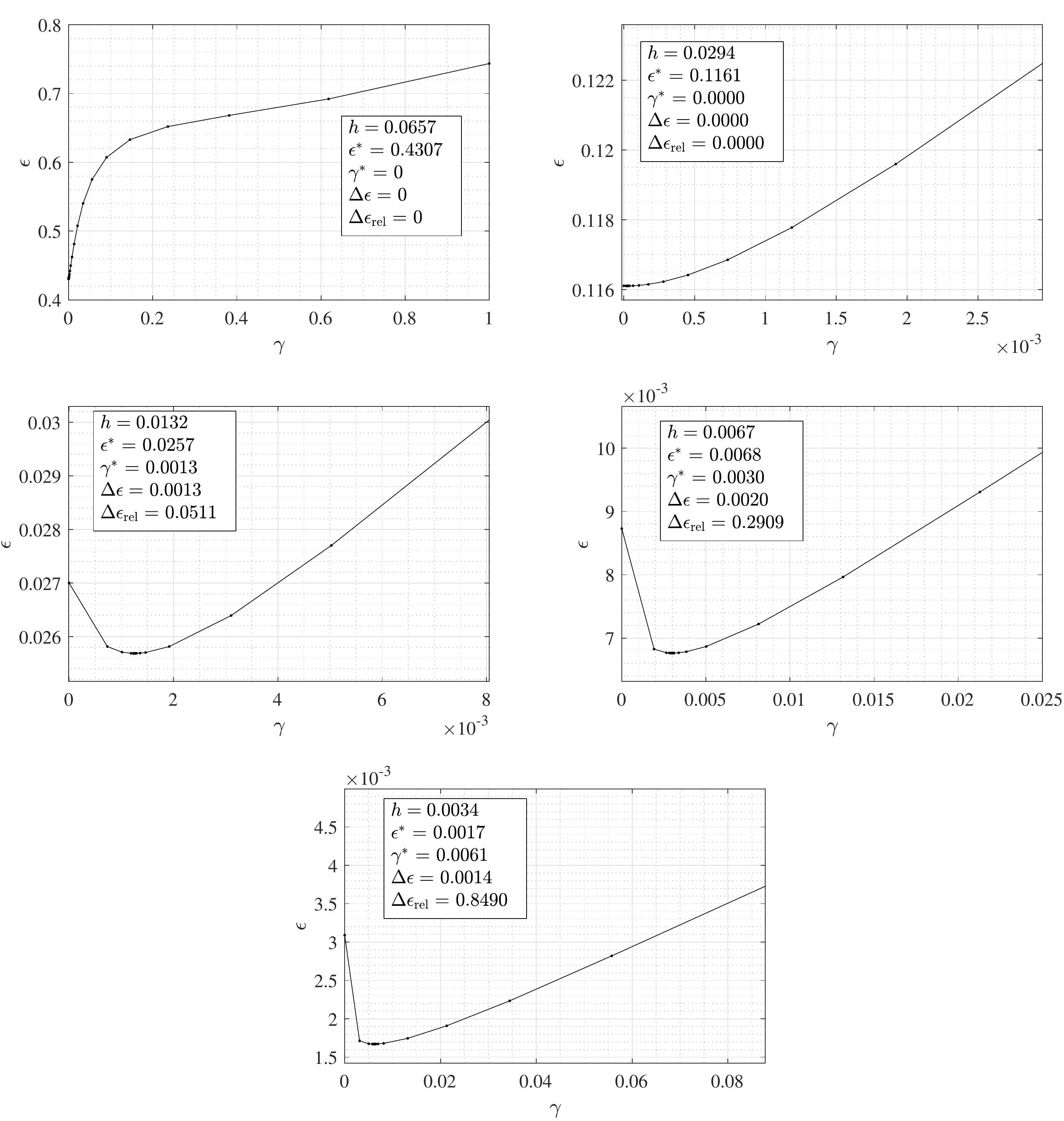}\caption{Vertex normal error as a function of the stabilization factor, $\epsilon(\gamma_{n}),$
for different mesh-sizes on a torus with $a=4$.\label{fig:gammaNormalA4}}
\end{figure}

In this section we analyze the influence of the stabilization factor
on the vertex normal error numerically by employing a golden search
method to find the optimal stabilization factor $\gamma_{n}^{*}$
that minimizes the normal error $\epsilon$ defined in (\ref{eq:VertexNormalErrorL2}).

\[
\begin{cases}
\underset{\gamma_{n}}{\mathrm{min}} & \epsilon(\gamma_{n})\\
\mathrm{s.t.} & \gamma_{n}^{0}\leq\gamma_{n}\leq\gamma_{n}^{1}
\end{cases},
\]
where we use $\gamma_{n}^{0}=0$ and $\gamma_{n}^{1}=1$. This is
done for several mesh-sizes and on a torus with $a=1$ and $a=4$,
see Figure \ref{fig:gammaNormalA1} and \ref{fig:gammaNormalA4}.
Notice how the curves become more planar, i.e., choosing a ``good''
$\gamma_{n}$ becomes less sensitive with the decrease in $h$.

\begin{figure}
\begin{centering}
\includegraphics[width=0.5\textwidth]{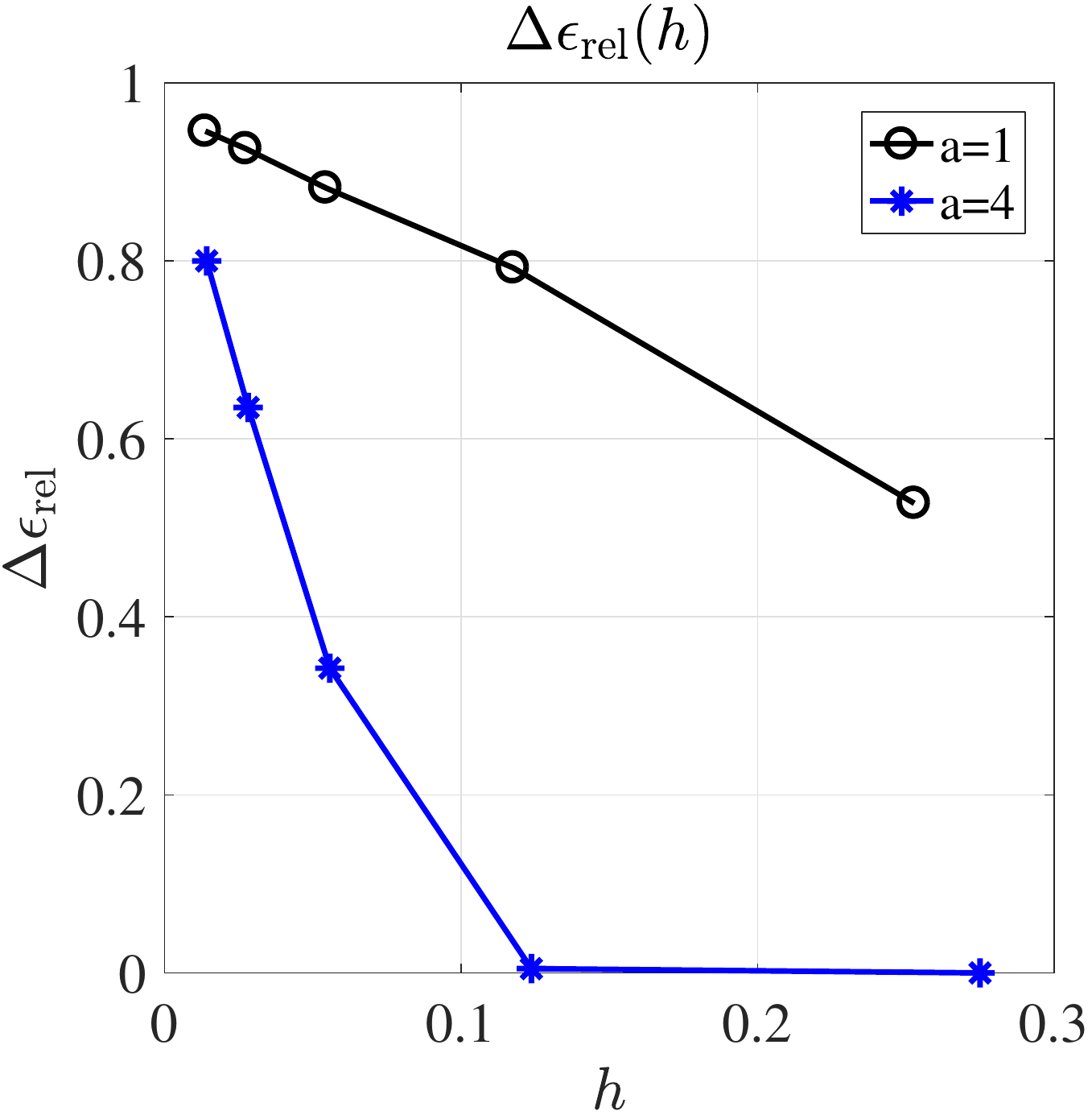}\caption{Relative difference in vertex normal error, $\Delta\epsilon_{\mathrm{rel}}$
as a function of the mesh-size $h$.\label{fig:DeltaEpsNormals}}
\par\end{centering}
\end{figure}

 The error difference is shown in Figure \ref{fig:DeltaEpsNormals}
and Table \ref{tab:gammaOfH} where $\epsilon_{L^{2}}$ is the $L^{2}$
error, defined in (\ref{eq:VertexNormalErrorL2}), of the $L^{2}$
vertex normals without stabilization and $\epsilon_{L_{\mathrm{stab}}^{2}}$
is the error of the stabilized $L^{2}$ vertex normal which is stabilized
with a optimal stabilization factor $\gamma^{*}$. 

\begin{table}
\begin{centering}
\subfloat[Torus with $a=1$]{\begin{centering}
\begin{tabular}{|c|c|c|c|}
\hline 
$h$ & $\gamma^{*}$ & $\epsilon_{L^{2}}-\epsilon_{L_{\mathrm{stab}}^{2}}|_{\gamma^{*}}$ & relative change\tabularnewline
\hline 
\hline 
0.0527 & 0.0147 & 0.0402 & 0.5277\tabularnewline
\hline 
0.0245 & 0.0371 & 0.0292 & 0.7923\tabularnewline
\hline 
0.0113 & 0.0693 & 0.0134 & 0.8821\tabularnewline
\hline 
0.0057 & 0.1057 & 0.0061 & 0.9262\tabularnewline
\hline 
0.0029 & 0.1623 & 0.0028 & 0.9459\tabularnewline
\hline 
\end{tabular}
\par\end{centering}
}
\par\end{centering}
\centering{}\subfloat[Torus with $a=4$]{\begin{centering}
\begin{tabular}{|c|c|c|c|}
\hline 
$h$ & $\gamma^{*}$ & $\epsilon_{L^{2}}-\epsilon_{L_{\mathrm{stab}}^{2}}|_{\gamma^{*}}$ & relative change\tabularnewline
\hline 
\hline 
0.0657 & 0 &  0 &  0\tabularnewline
\hline 
0.0294 & 0.0000 &     0.0003 &     0.0048\tabularnewline
\hline 
0.0132 & 0.0013 &     0.0078 &     0.3424\tabularnewline
\hline 
0.0067 & 0.0030 &     0.0065 &     0.6354\tabularnewline
\hline 
0.0034 & 0.0061 &     0.0035 &     0.8001\tabularnewline
\hline 
\end{tabular}
\par\end{centering}
}\caption{Stabilization factor as a function of mesh-size.\label{tab:gammaOfH} }
\end{table}

\subsection{Interpolation}

\begin{figure}
\begin{centering}
\includegraphics[width=0.8\textwidth]{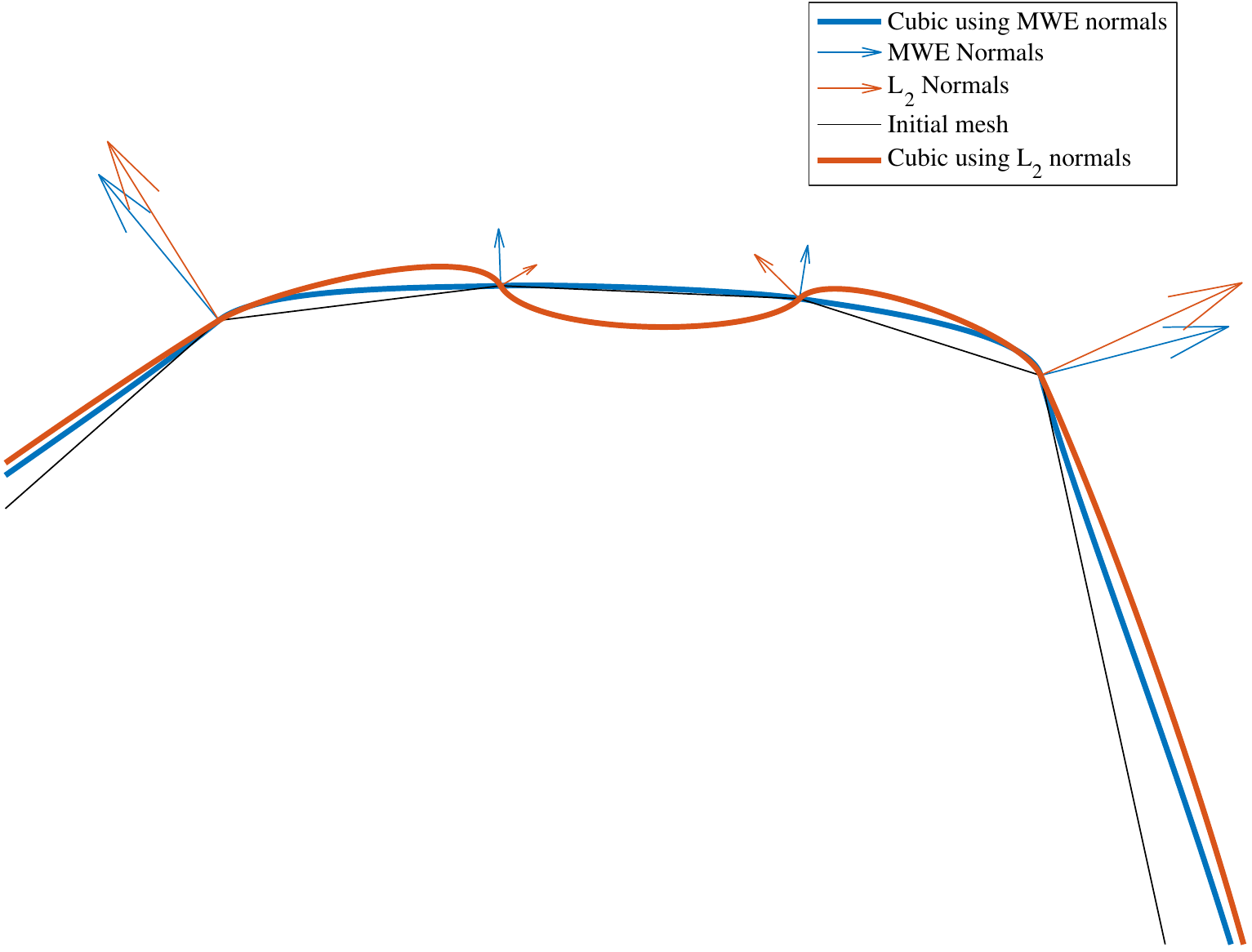}\caption{2D cubic Hermite interpolation of a coarse line segment using two
different approximations of the vertex normals. \label{fig:2D-cubic-Hermite}}
\par\end{centering}
\end{figure}

In a 2D case we can see in Figure \ref{fig:2D-cubic-Hermite} how
the choice of vertex normals affects the resulting cubic interpolation.
The initial mesh is coarse and the (unstabilized) $L^2$--projected normals
are not just depending on the nearest neighbors to each vertex but
globally. The resulting difference is apparent.

\begin{figure}
\begin{centering}
\includegraphics[width=0.8\textwidth]{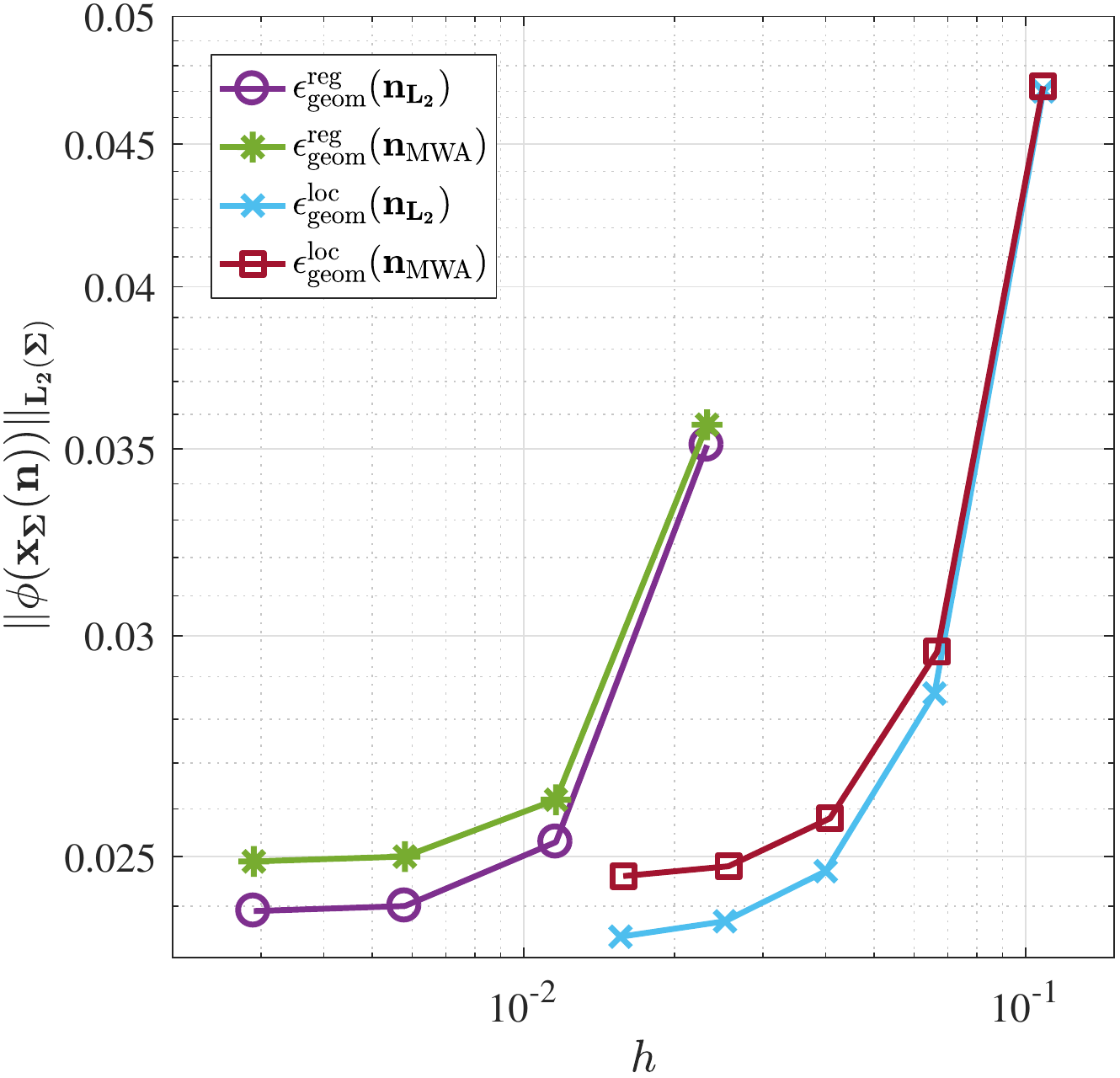}\caption{Geometrical error of a torus with $a=4$, initial mesh-size of 0.1618
and initial geometrical error of 0.0863 \label{fig:GeoErrorPNRef}}
\par\end{centering}
\end{figure}

We compare the impact of different vertex normals on the interpolation
by measuring the geometrical error defined as

\[
\epsilon_{\mathrm{geom}}=\|\phi(\bm{x}_{\Sigma}(\bm{n}))^{2}\|_{L^2_h},
\]
where $\bm{x}_{\Sigma}(\bm{n})$ denotes the discrete surface interpolated
with a particular normal approximation method. We measure the $L^2_h$--norm of the signed distance. The refinement algorithm employed is
the PN triangles using 1 tessellation per face see Figure \ref{fig:PNTriangleN1}.
The mesh-size in this section is defined as 
\[
h:=\left(\sum_{K=1}^{N_{\text{e}}}\sqrt{A_{K}}\right)/N_{\text{e}},
\]
where $N_{\text{e}}$ denotes the number of elements and $A_{K}$ is the
area of the $K$-th element. The initial mesh-size is $h=0.1618$
and the initial $L^2_h$--norm of the signed distance error is $\epsilon_{\mathrm{geom}}=0.0863$.
See Figure \ref{fig:GeoErrorPNRef} for the convergence comparison,
Table \ref{tab:GeoErrorPNRef} for the regular refinement data and
Table \ref{tab:GeoErrorPNLocRef} for the local refinement data. 

\begin{figure}
\begin{centering}
\subfloat[]{\begin{centering}
\includegraphics[width=0.49\textwidth]{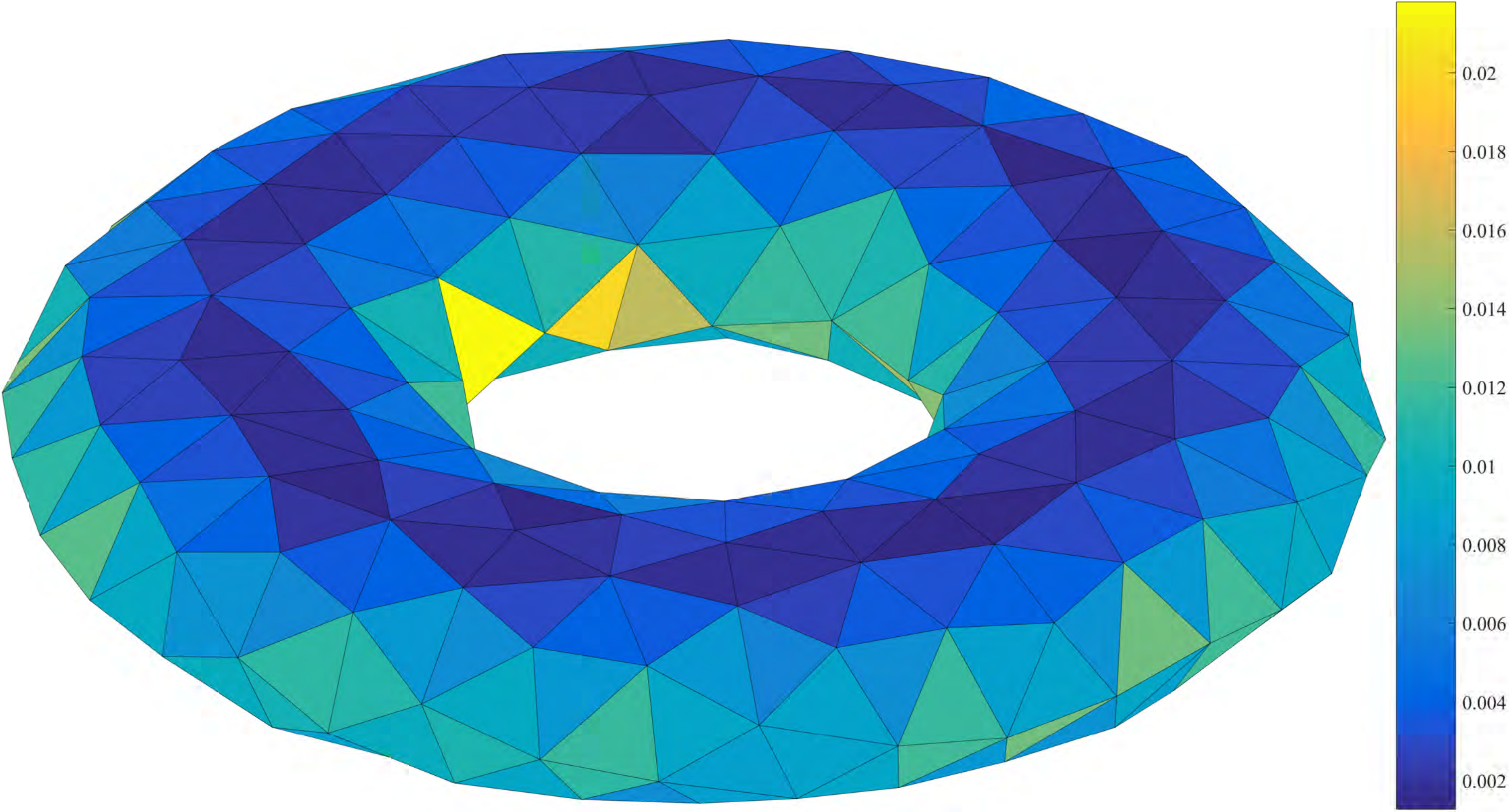}
\par\end{centering}

}\subfloat[]{\begin{centering}
\includegraphics[width=0.49\textwidth]{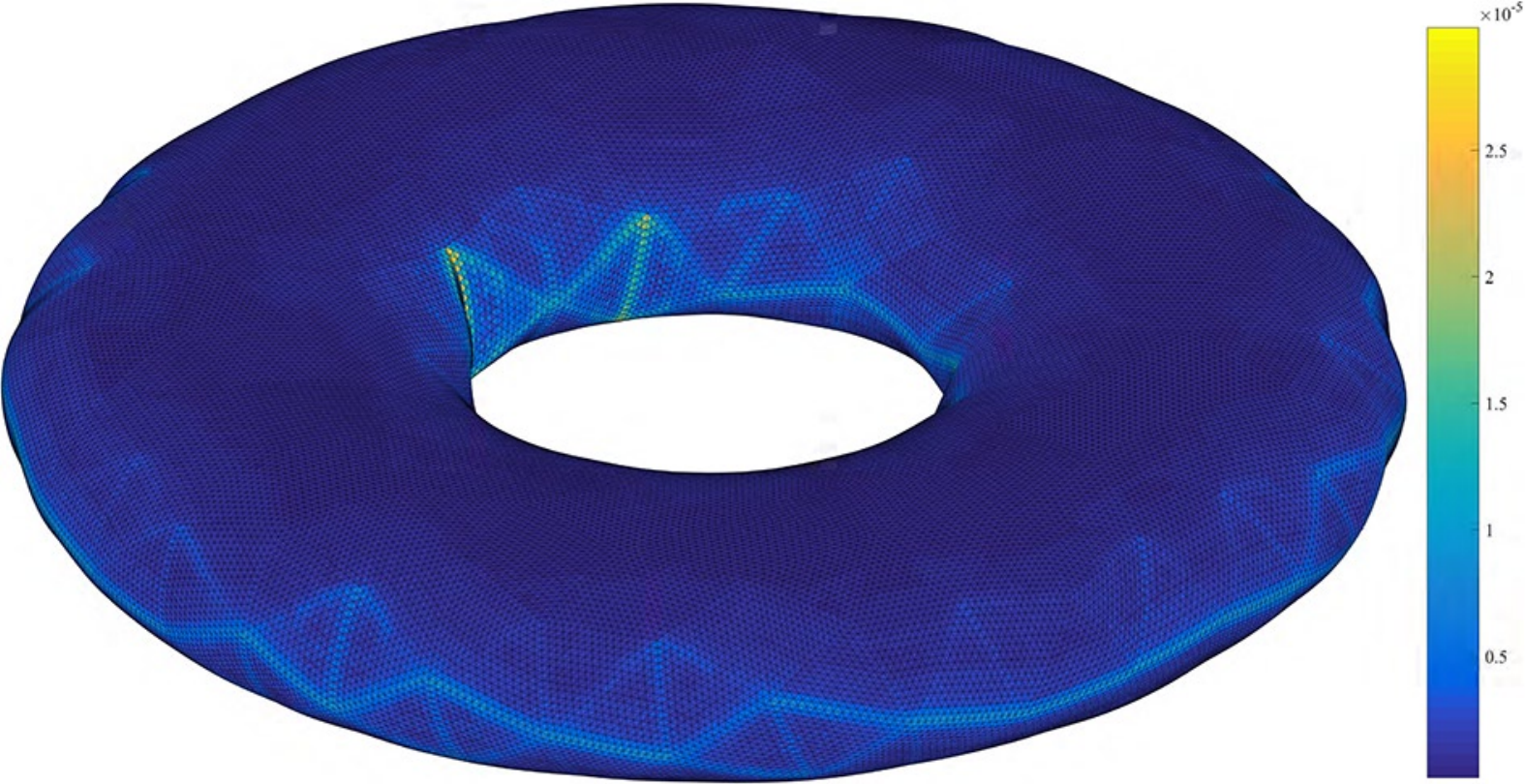}
\par\end{centering}
}
\par\end{centering}
\begin{centering}
\subfloat[]{\begin{centering}
\includegraphics[width=0.49\textwidth]{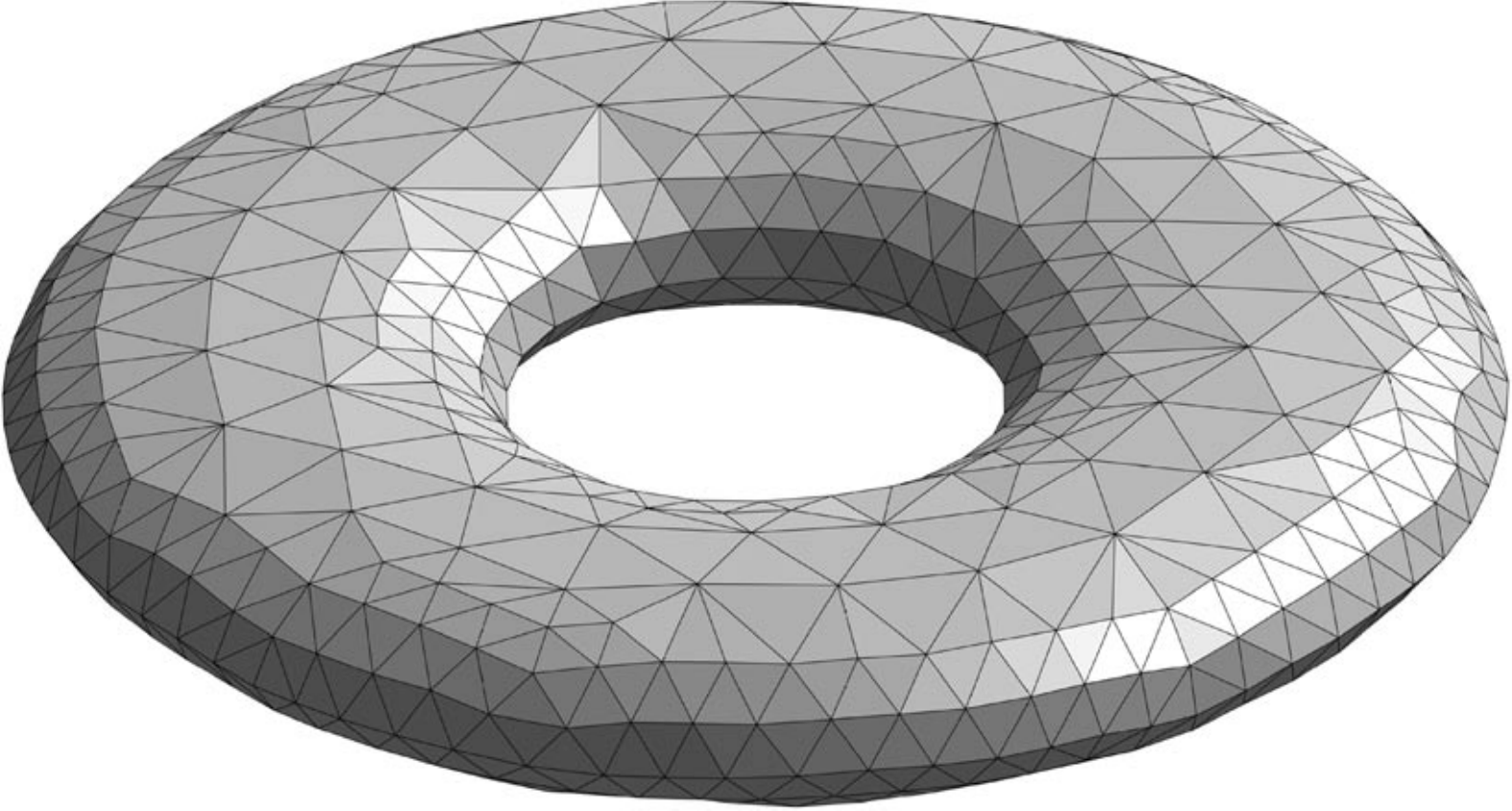}
\par\end{centering}
}\subfloat[]{\begin{centering}
\includegraphics[width=0.49\textwidth]{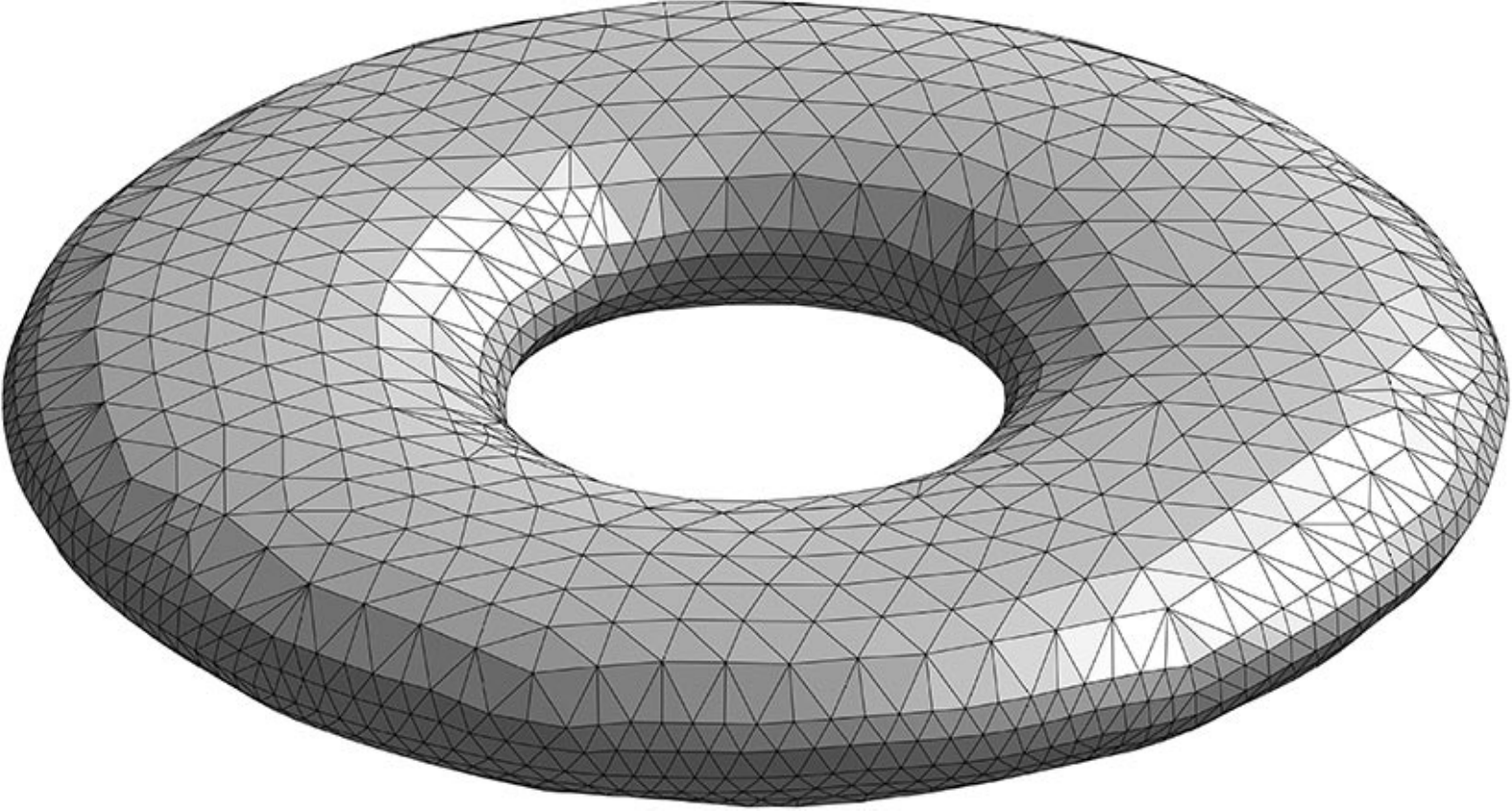}
\par\end{centering}
}
\par\end{centering}
\begin{centering}
\subfloat[]{\begin{centering}
\includegraphics[width=0.49\textwidth]{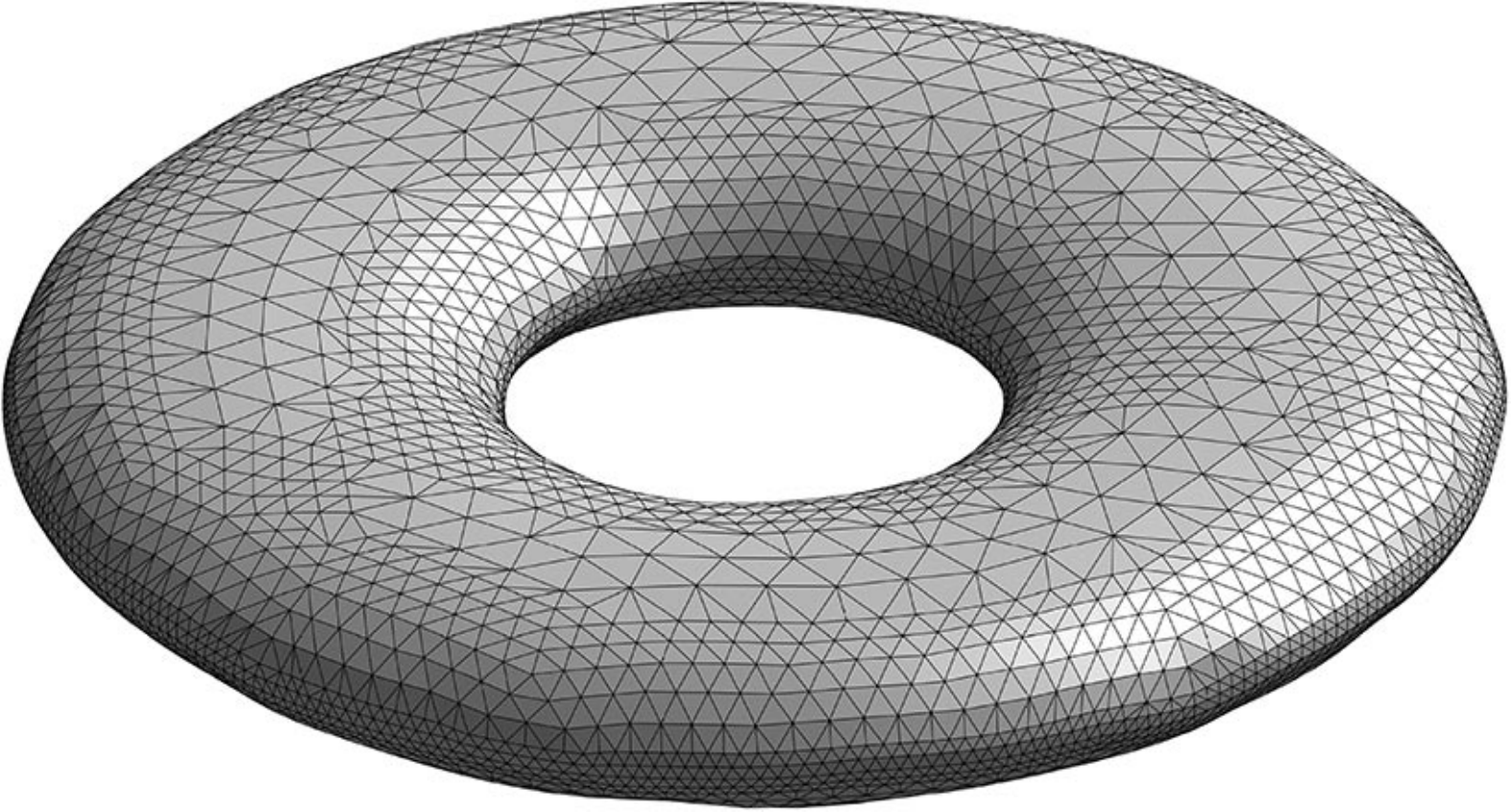}
\par\end{centering}
}\subfloat[]{\begin{centering}
\includegraphics[width=0.49\textwidth]{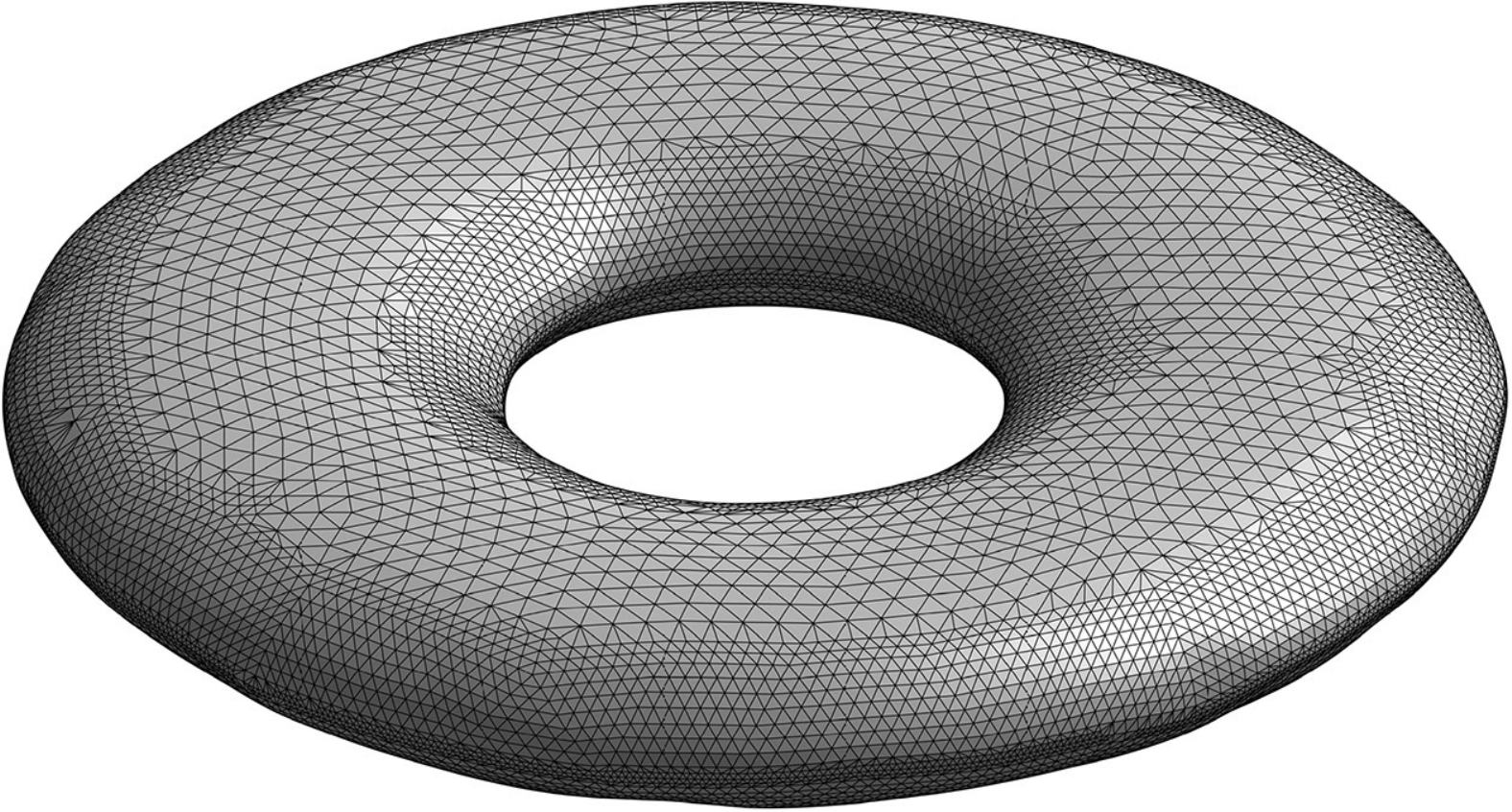}
\par\end{centering}
}\caption{PN refinement on the torus using $a=4$ and $N=1$. The plot is showing
$\|\bm{n}_{f}-\bm{n}_{{L^2_{\mathrm{stab}}}}\|_{L^2_{h}}$
on: (a) Initial unrefined torus. (b) Torus after 4 regular refinements.
(c)-(f) Torus local refinements 1 to 4.\label{fig:PN-Interpolation-local-refinement-Torus}}
\par\end{centering}
\end{figure}

\begin{figure}
\begin{centering}
\subfloat[]{
\begin{centering}
\includegraphics[width=0.49\textwidth]{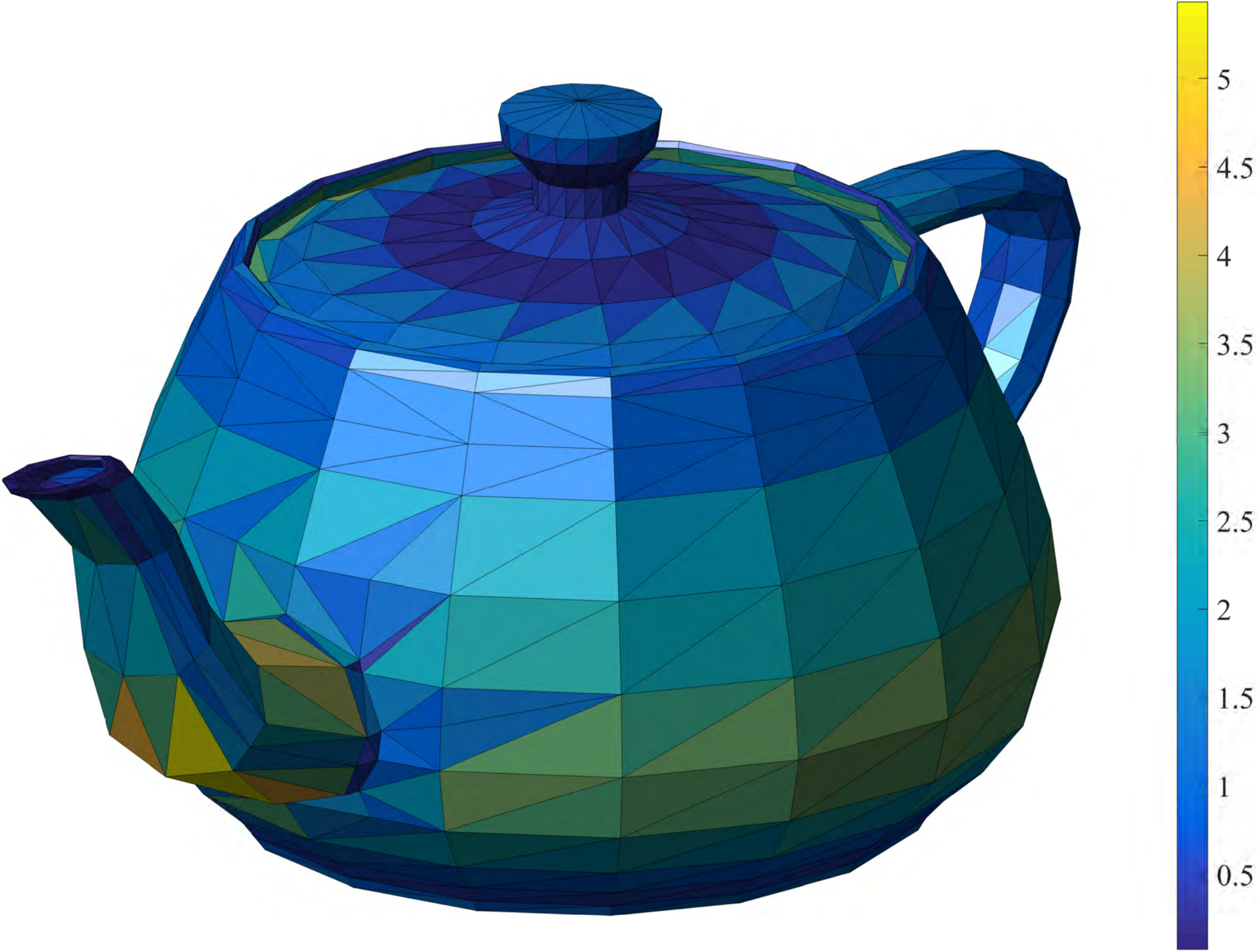}
\par\end{centering}
}\subfloat[]{\centering{}\includegraphics[width=0.49\textwidth]{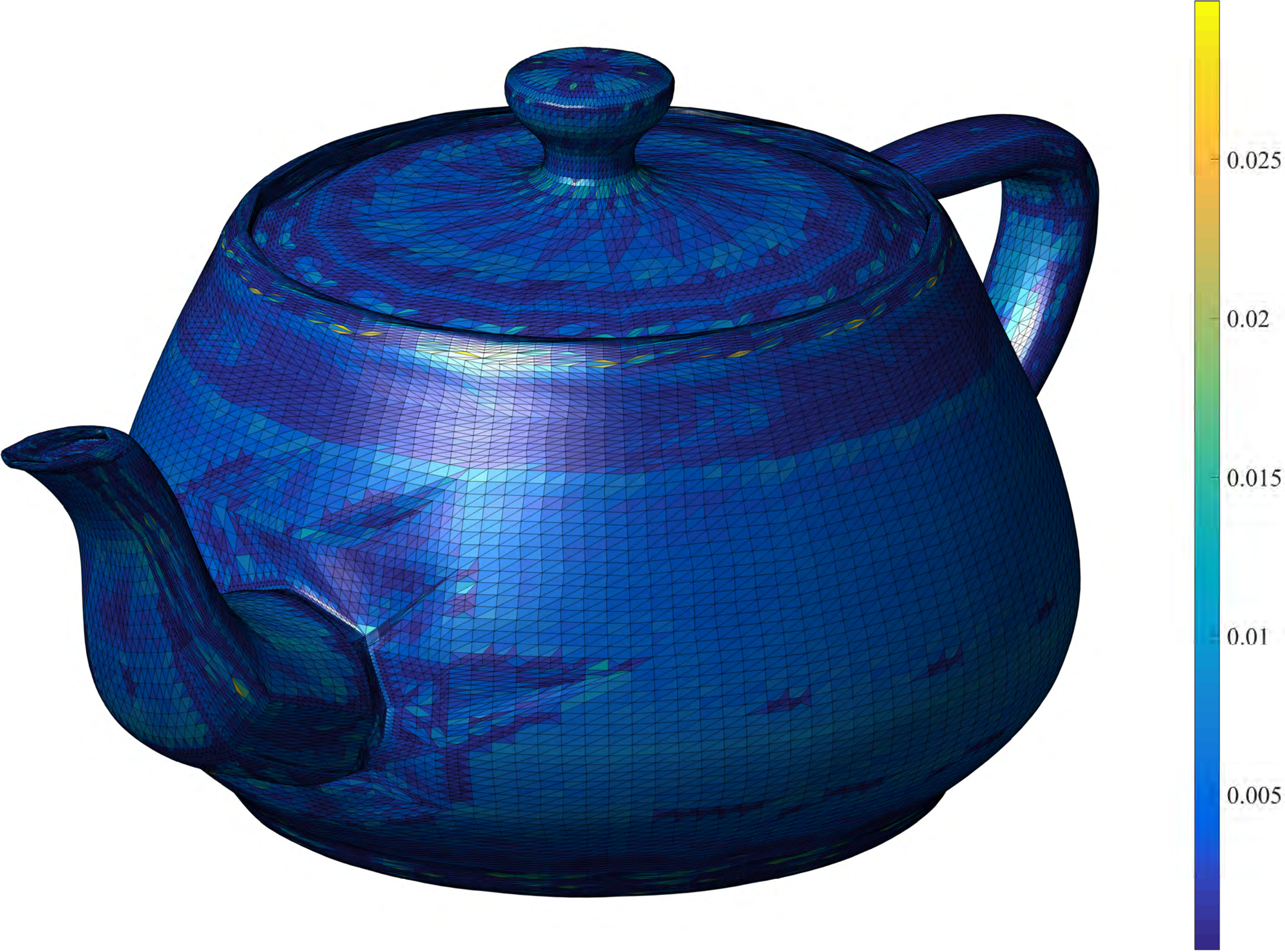}}\caption{Closed surface of a teapot mesh of fairly low quality. The plot is
showing $\|\bm{n}_{f}-\bm{n}_{{L^2_{\mathrm{stab}}}}\|_{L^2_{h}}$
on: (a) Initial surface. (b) After 4 local refinements. \label{fig:PN-local-refinement Teapot}}
\par\end{centering}
\end{figure}

\begin{figure}
\centering{}\subfloat[]{\centering{}\includegraphics[width=0.49\textwidth]{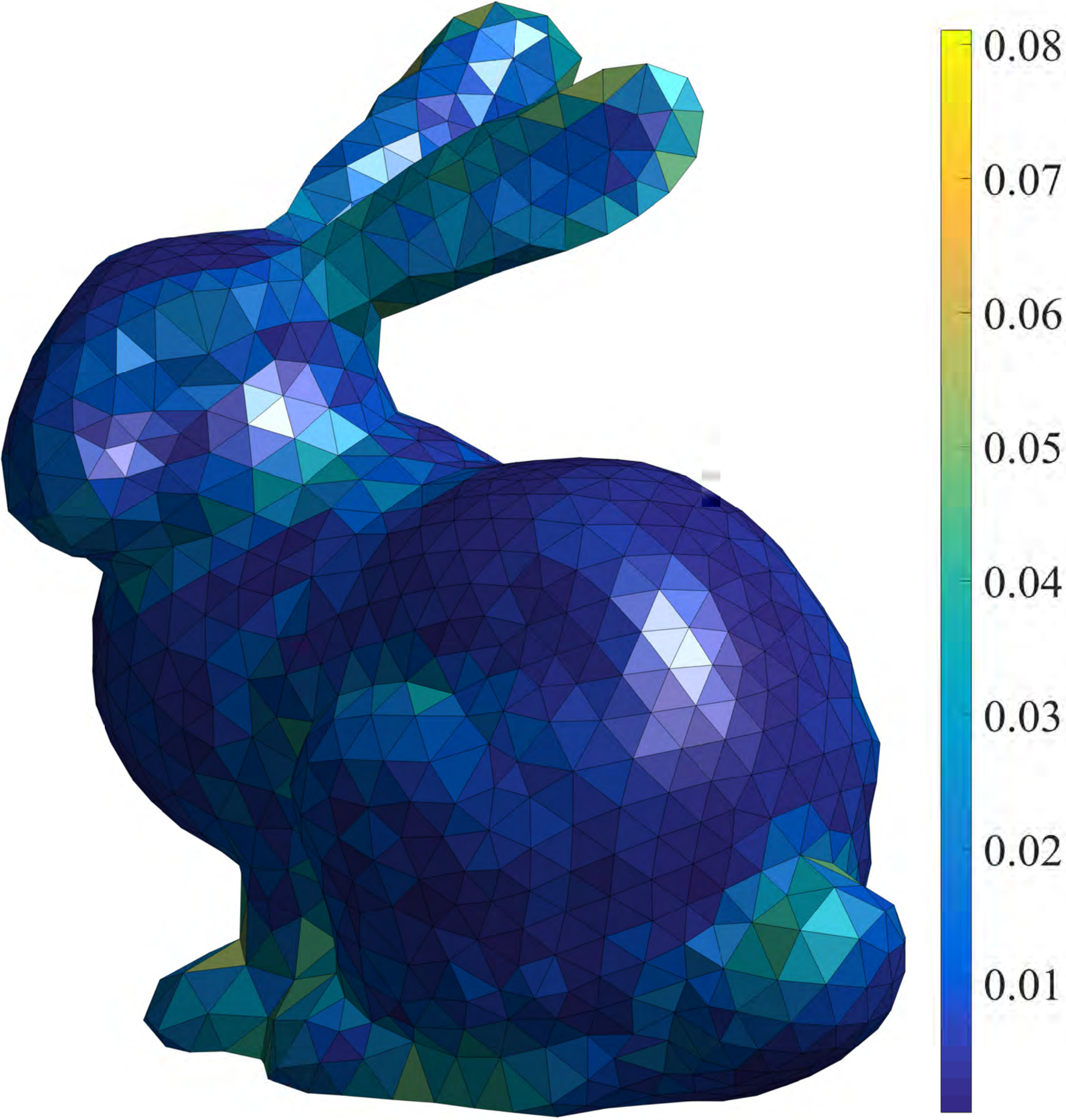}}\subfloat[]{\centering{}\includegraphics[width=0.49\textwidth]{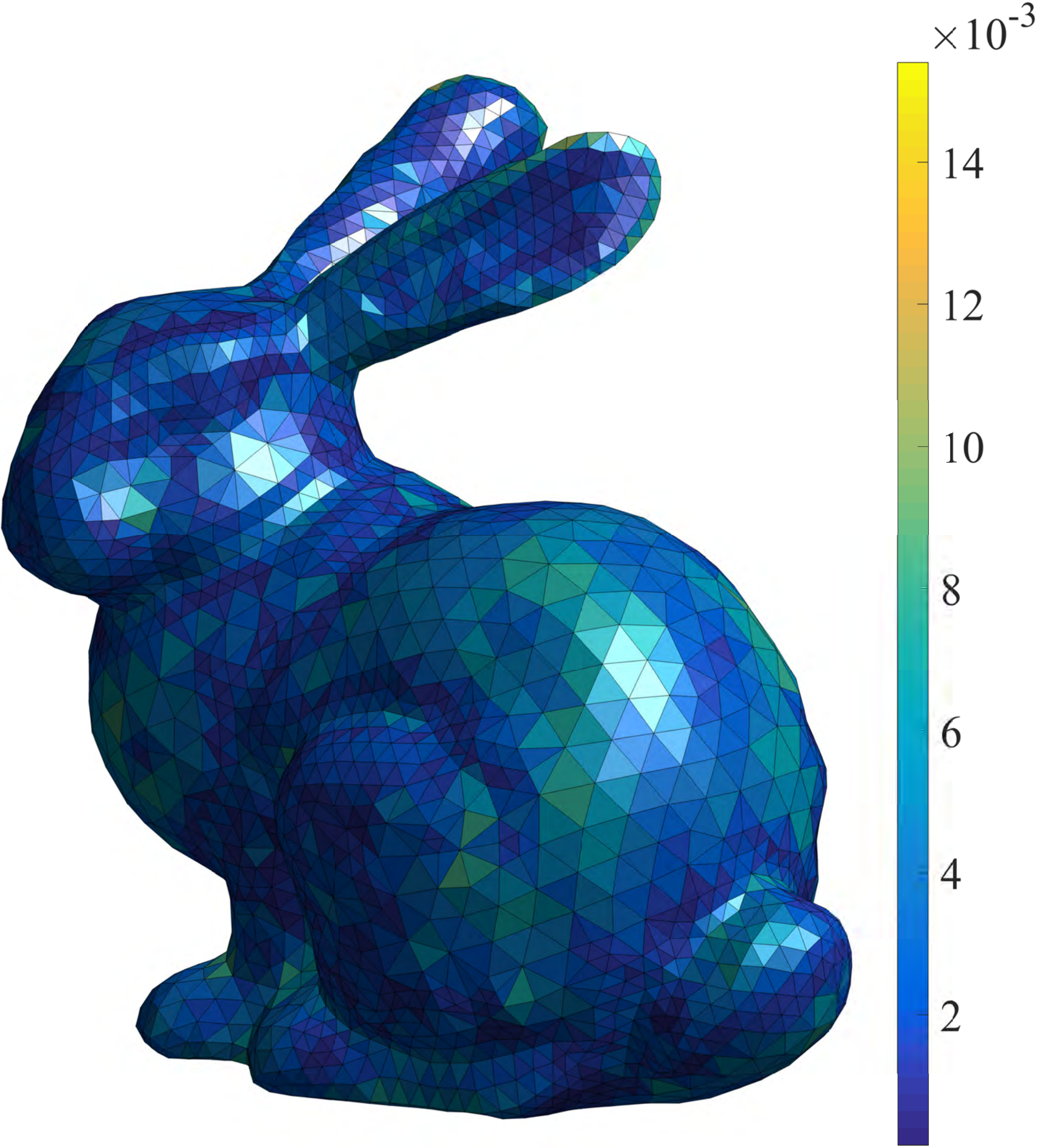}}\caption{Stanford bunny. The plot is showing $\|\bm{n}_{f}-\bm{n}_{{L^2_{\mathrm{stab}}}}\|_{L^2_{h}}$
on: (a) Initial surface. (b) After 1 local refinement.  \label{fig:PN-local-refinement Bunny}}
\end{figure}

Examples of interpolation using PN triangles with local refinement are shown in Figure \ref{fig:PN-Interpolation-local-refinement-Torus} for the Torus, Figure \ref{fig:PN-local-refinement Teapot} for the Utah teapot and Figure \ref{fig:PN-local-refinement Bunny} for the Stanford bunny.

\begin{figure}
\begin{centering}
\subfloat[]{\begin{centering}
\includegraphics[width=0.49\textwidth]{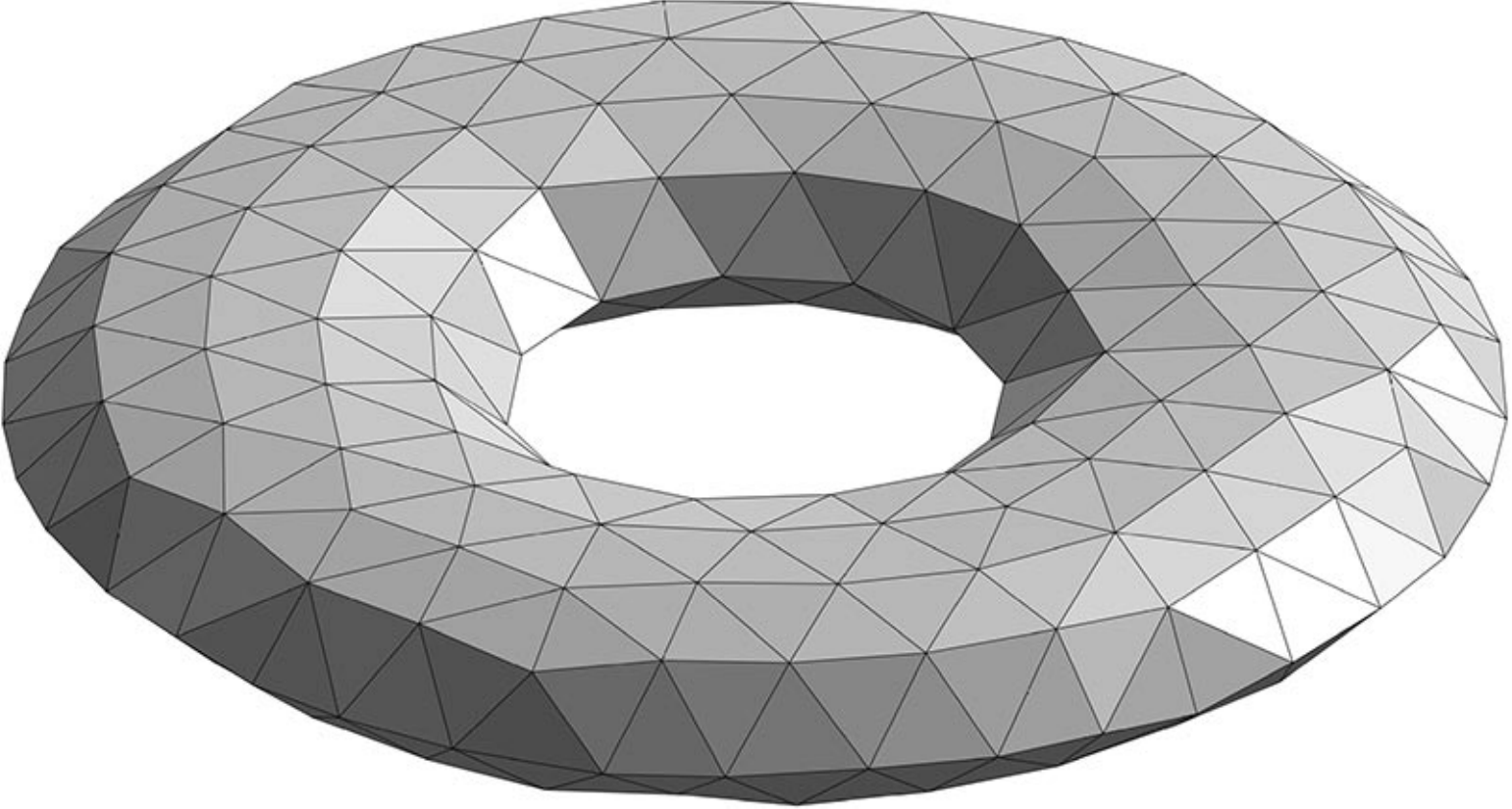}
\par\end{centering}
}\subfloat[]{\begin{centering}
\includegraphics[width=0.49\textwidth]{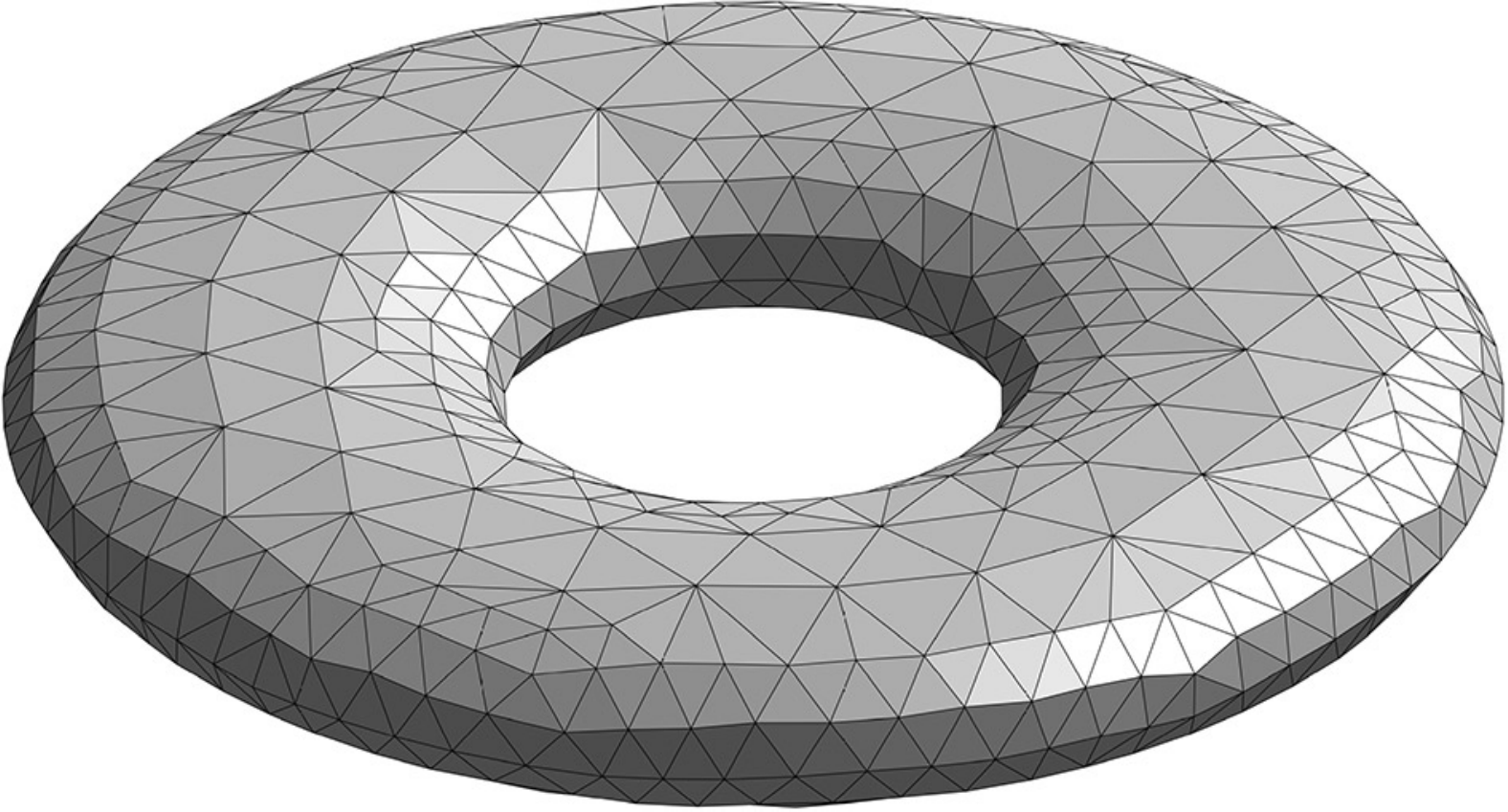}
\par\end{centering}
}
\par\end{centering}
\begin{centering}
\subfloat[]{\begin{centering}
\includegraphics[width=0.49\textwidth]{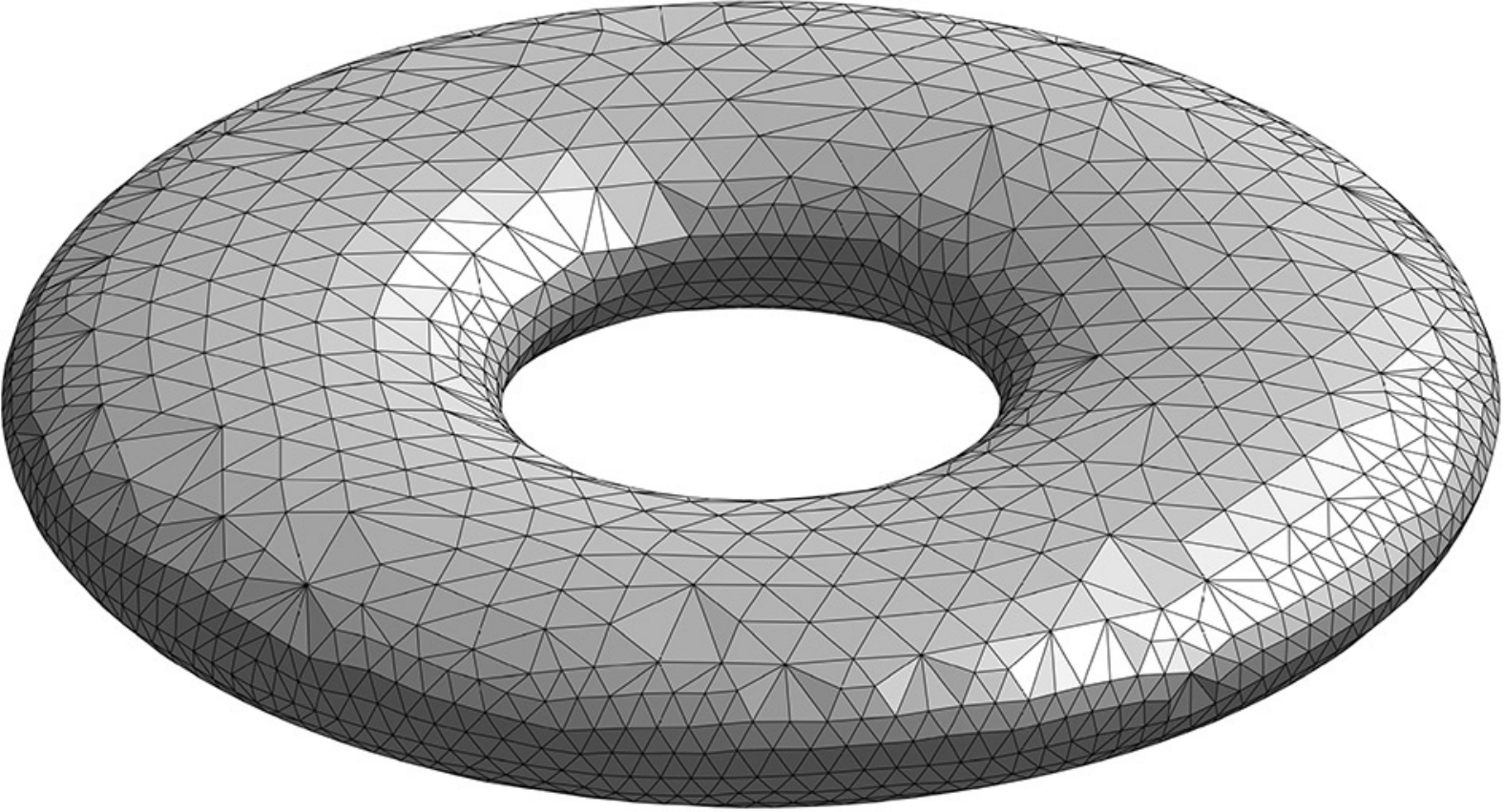}
\par\end{centering}
}\subfloat[]{\begin{centering}
\includegraphics[width=0.49\textwidth]{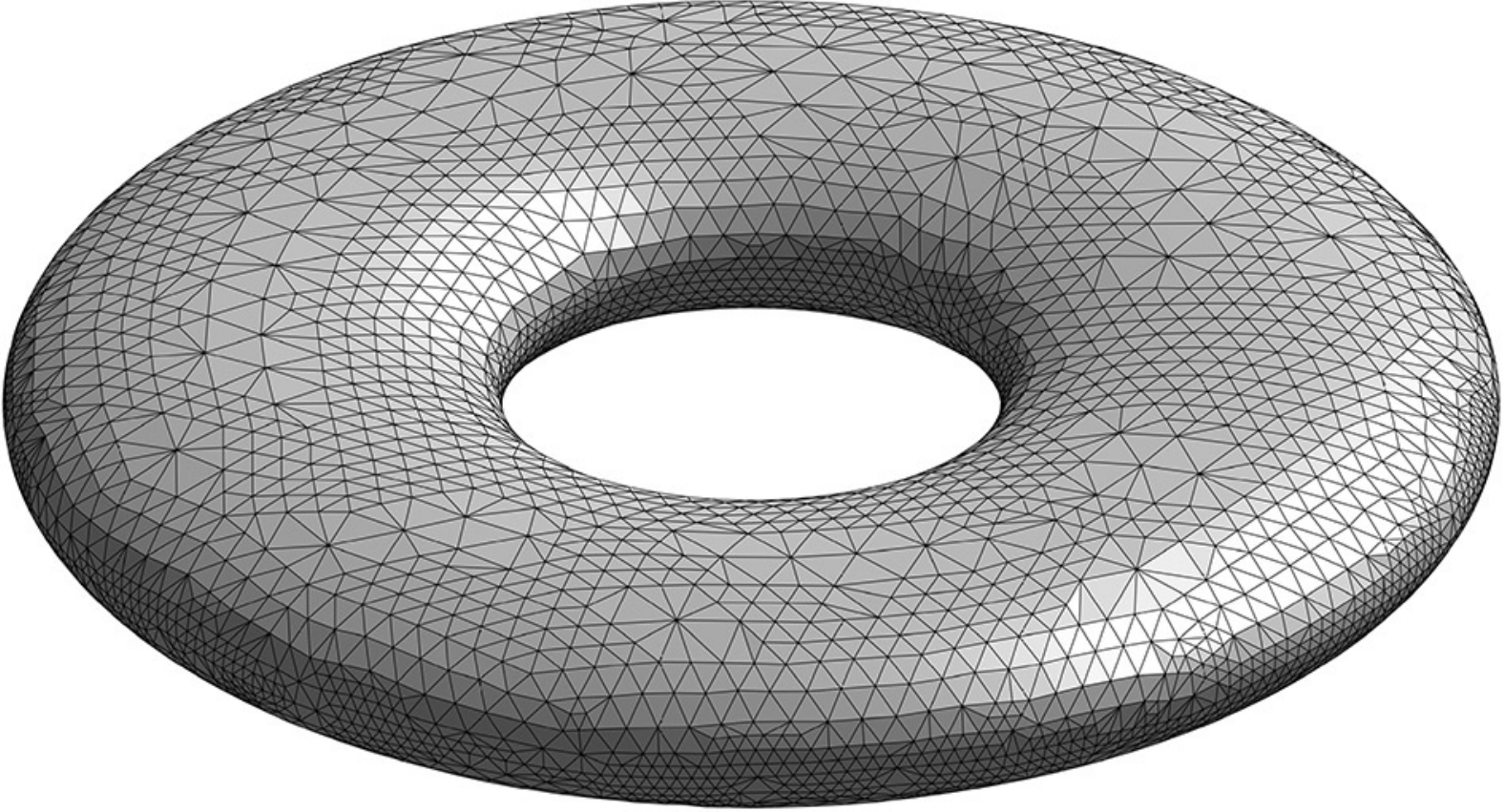}
\par\end{centering}
}
\par\end{centering}
\centering{}\subfloat[]{\begin{centering}
\includegraphics[width=0.49\textwidth]{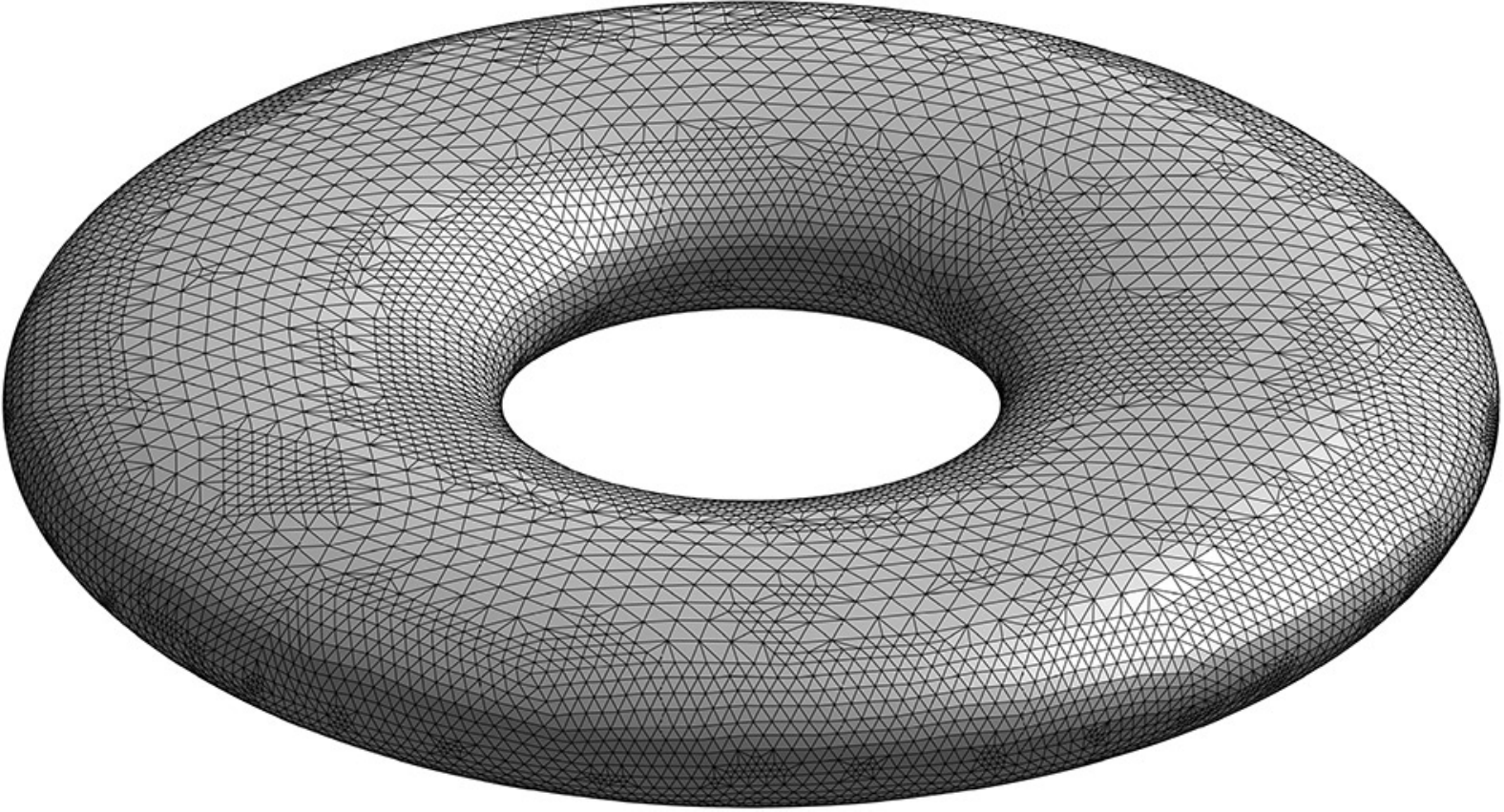}
\par\end{centering}
}\caption{PN local refinement with projection to the exact surface on the torus
using $a=4$ and $N=1$. \label{fig:PN-local-refinement Torus-Exact}}
\end{figure}

The local refinement method is compared to local refinement with projection
to the exact surface, see Figure \ref{fig:PN-local-refinement Torus-Exact}.
We compare the approximite normal error with the exact normal error
by computing the effectivity index, given by

\[
E=\dfrac{\|\bm{n}_{e}-\bm{n}_{f}\|_{L^2_{h}}}{\|\bm{n}_{L^2_{\mathrm{stab}}}-\bm{n}_{f}\|_{L^2_{h}}}
\]
where $\bm{n}_{e}$ is the exact normal to the surface, $\bm{n}_{f}$
is the face normal and $\bm{n}_{L^2_{\mathrm{stab}}}$ is the recovered
stabilized $L^2$--projected normal, see Table \ref{tab:GeoErrorPNLocRefExact}.

\begin{table}
\centering{} \resizebox{\textwidth}{!}{
\begin{tabular}{|c|c|c|c|c|c|}
\hline 
$N_{e}$ & $h$ & $\epsilon_{\mathrm{geom}}|_{\bm{n}_{\mathrm{L2stab}}}$ & $\|\bm{n}_{f}-\bm{n}_{\textrm{\ensuremath{\mathrm{L2stab}}}}\|_{\Sigma_{h}}$ & $\epsilon_{\mathrm{geom}}|_{\bm{n}_{\mathrm{MWA}}}$ & $\|\bm{n}_{f}-\bm{n}_{\textrm{\ensuremath{\mathrm{MWA}}}}\|_{\Sigma_{h}}$\tabularnewline
\hline 
\hline 
1856 & 0.0900 & 0.0351 & 1.3934 & 0.0357 & 1.4351\tabularnewline
\hline 
7424 & 0.0451 & 0.0253 & 0.6976 & 0.0262 & 0.6981\tabularnewline
\hline 
29696 & 0.0226 & 0.0240 & 0.3453 & 0.0250 & 0.3354\tabularnewline
\hline 
118784 & 0.0113 & 0.0239 & 0.1766 & 0.0249 & 0.1619\tabularnewline
\hline 
\end{tabular}}\caption{PN regular refinement of a torus with $a=4$, initial mesh-size of
0.0321 and initial geometrical error of 0.0863.\label{tab:GeoErrorPNRef} }
\end{table}

\begin{table}
\begin{centering}
 \resizebox{\textwidth}{!}{
\begin{tabular}{|c|c|c|c|c|c|c|c|}
\hline 
$N_{e}$ & $h$ & $\epsilon_{\mathrm{geom}}|_{\bm{n}_{\mathrm{L2stab}}}$ & $\|\bm{n}_{f}-\bm{n}_{\textrm{\ensuremath{\mathrm{L2stab}}}}\|_{\Sigma_{h}}$ & $N_{e}$ & $h$ & $\epsilon_{\mathrm{geom}}|_{\bm{n}_{\mathrm{MWA}}}$ & $\|\bm{n}_{f}-\bm{n}_{\textrm{\ensuremath{\mathrm{MWA}}}}\|_{\Sigma_{h}}$\tabularnewline
\hline 
\hline 
1278 & 0.1085 & 0.0470 & 1.6393 & 1278 & 0.1085 & 0.0472 & 1.6831\tabularnewline
\hline 
3470 & 0.0659 & 0.0286 & 0.9966 & 3390 & 0.0667 & 0.0296 & 1.0123\tabularnewline
\hline 
9420 & 0.0400 & 0.0247 & 0.5945 & 9068 & 0.0408 & 0.0258 & 0.5952\tabularnewline
\hline 
23862 & 0.0252 & 0.0237 & 0.3675 & 23060 & 0.0256 & 0.0248 & 0.3643\tabularnewline
\hline 
62068 & 0.0156 & 0.0234 & 0.2285 & 60294 & 0.0158 & 0.0246 & 0.2205\tabularnewline
\hline 
\end{tabular}}\caption{Local refinement of a torus with $a=4$, initial mesh-size of 0.0321
and initial geometrical error of 0.0863.\label{tab:GeoErrorPNLocRef}}
\par\end{centering}
\end{table}

\begin{table}
\centering{} \resizebox{\textwidth}{!}{
\begin{tabular}{|c|c|c|c|c|c|}
\hline 
$N_{e}$ & $h$ & $\epsilon_{\mathrm{geom}}|_{\bm{n}_{\mathrm{L2stab}}}$ & Rate & $\|\bm{n}_{f}-\bm{n}_{\textrm{\ensuremath{\mathrm{L2stab}}}}\|_{\Sigma_{h}}$ & Effectivity Index, $E$\tabularnewline
\hline 
\hline 
1278 & 0.0635 & 0.0393 & - & 1.5308 & 1.0055\tabularnewline
\hline 
3474 & 0.0307 & 0.0143 & 1.3900 & 0.8872 & 1.0147\tabularnewline
\hline 
8736 & 0.0151 & 0.0054 & 1.3608 & 0.5671 & 1.0024\tabularnewline
\hline 
21149 & 0.0095 & 0.0025 & 1.7165 & 0.3527 & 1.0013\tabularnewline
\hline 
53047 & 0.0056 & 0.0010 & 1.7790 & 0.2225 & 1.0009\tabularnewline
\hline 
131740 & 0.0028 & 0.0004 & 1.3741 & 0.1429 & 0.9999\tabularnewline
\hline 
\end{tabular}}\caption{Local refinement with projection to the exact surface of a torus with
$a=4$, initial mesh-size of 0.1618.\label{tab:GeoErrorPNLocRefExact}}
\end{table}

\subsection{Mean curvature}

The mean curvature is computed on a structured and unstructured torus
with $R=1$, $r=1/2$ and $a=1$. We compare the mean curvature approximation
to the exact mean curvature, the smooth surface fit approach (SSF)
and the discrete local Laplace-Beltrami (DLLB) approach described
in Section \ref{subsec:DiscreteLocalLaplaceBeltrami}, and our stabilized discrete curvature vector solving (\ref{stabilizedcurvature}).
In the last case we compute the mean curvature $H_{h}$ through
\[
H_{h}=\frac{\bm{H}_{h}\cdot\bm{n}_{h}}{2},
\]
where $\bm{n}_{h}$ denotes the normal computed using the
stabilized $L^{2}$\textendash projection from (\ref{eq:discreteStabVertexNormal}).
In our computational experience, this gives a more accurate result
than the immediate $H_{h}=\frac{1}{2}\vert\bm{H}_{h}\vert$.

\begin{figure}
\begin{centering}
\subfloat[Exact]{\centering{}\includegraphics[width=0.5\columnwidth]{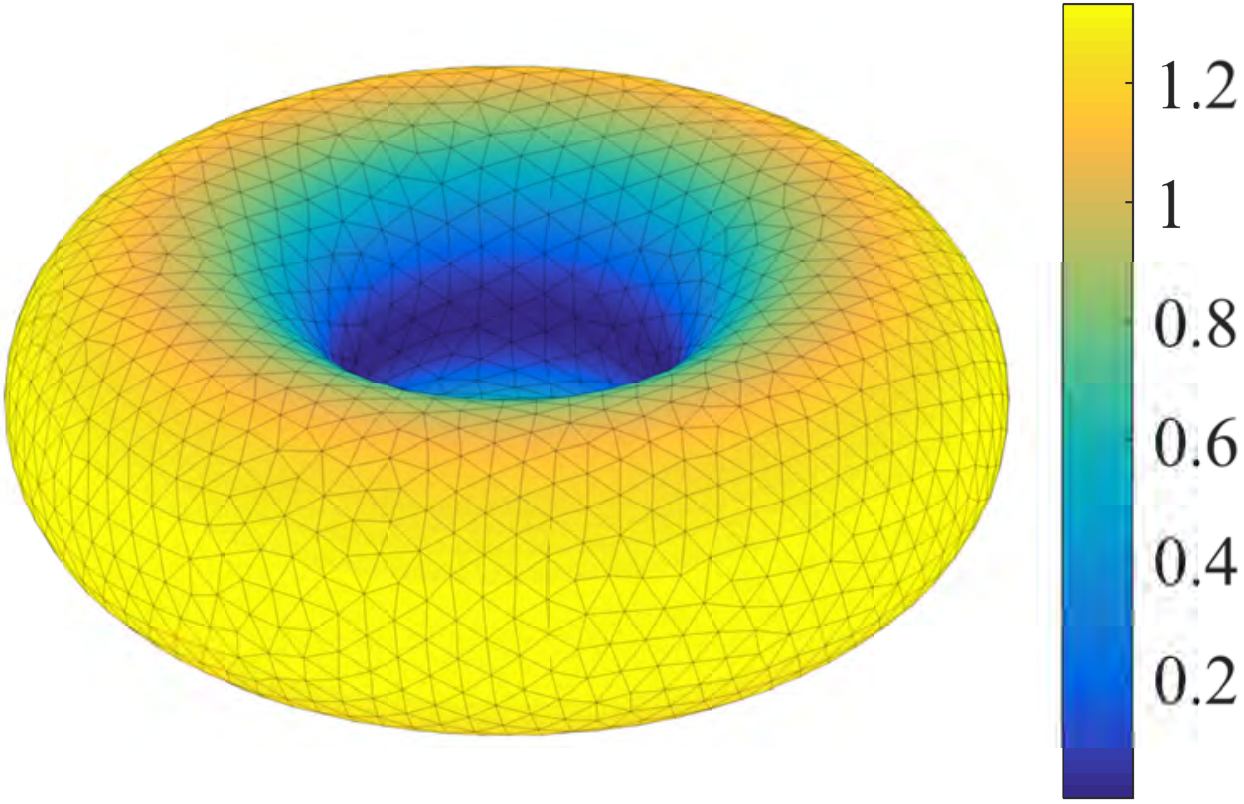}}\subfloat[SSF]{\centering{}\includegraphics[width=0.5\columnwidth]{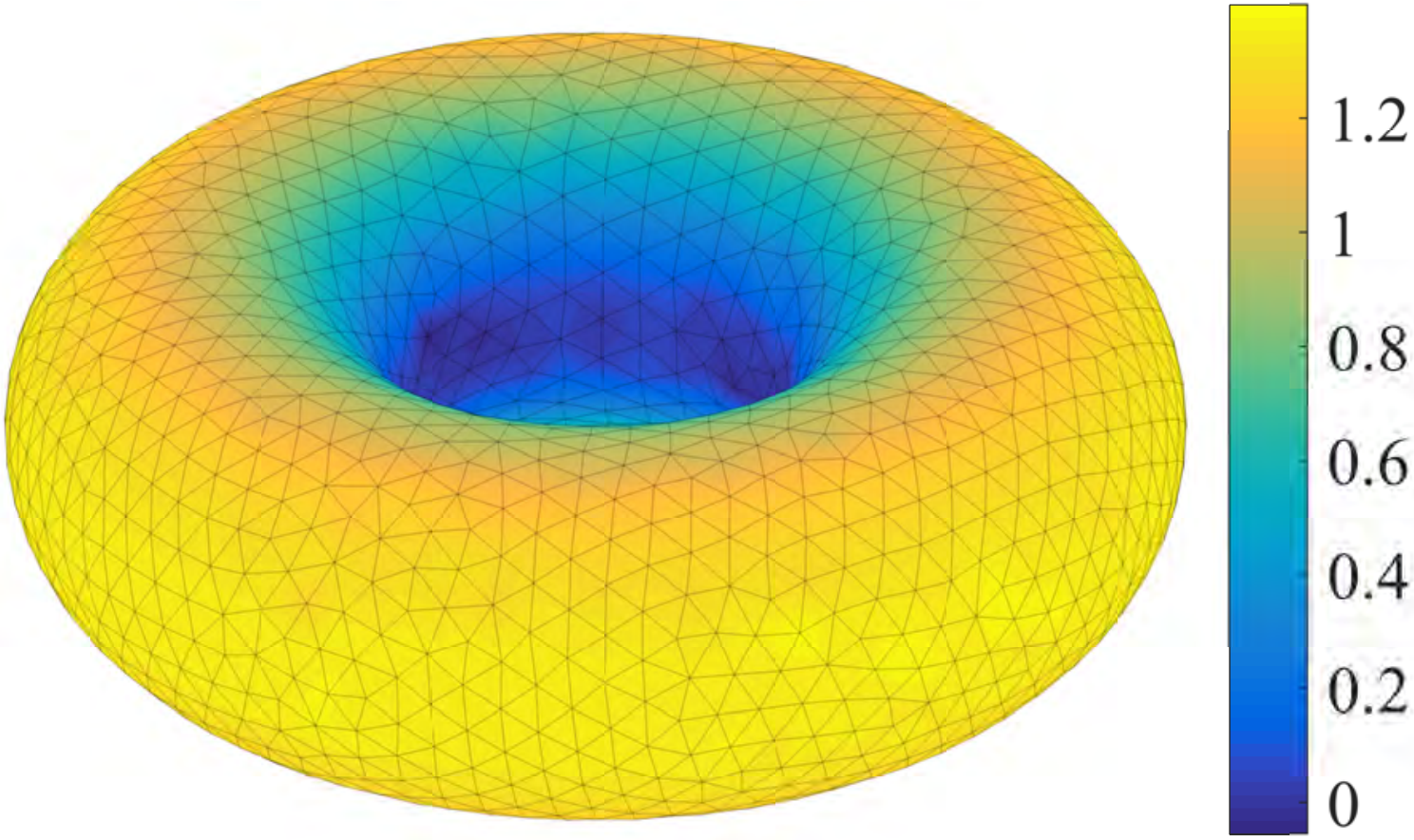}}
\par\end{centering}
\begin{centering}
\subfloat[DLLB]{\centering{}\includegraphics[width=0.5\columnwidth]{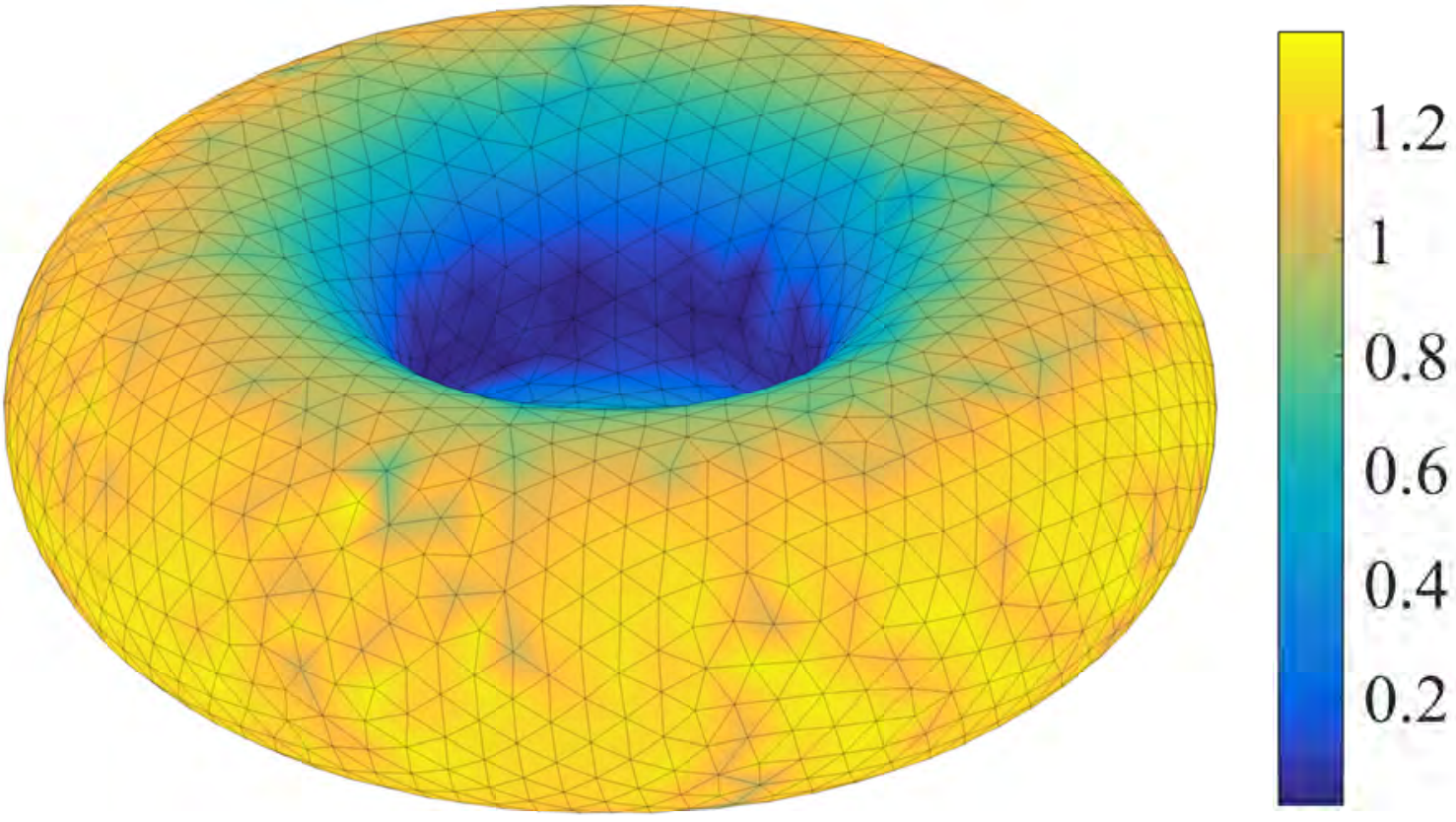}}\subfloat[$H_{h}$]{\centering{}\includegraphics[width=0.5\columnwidth]{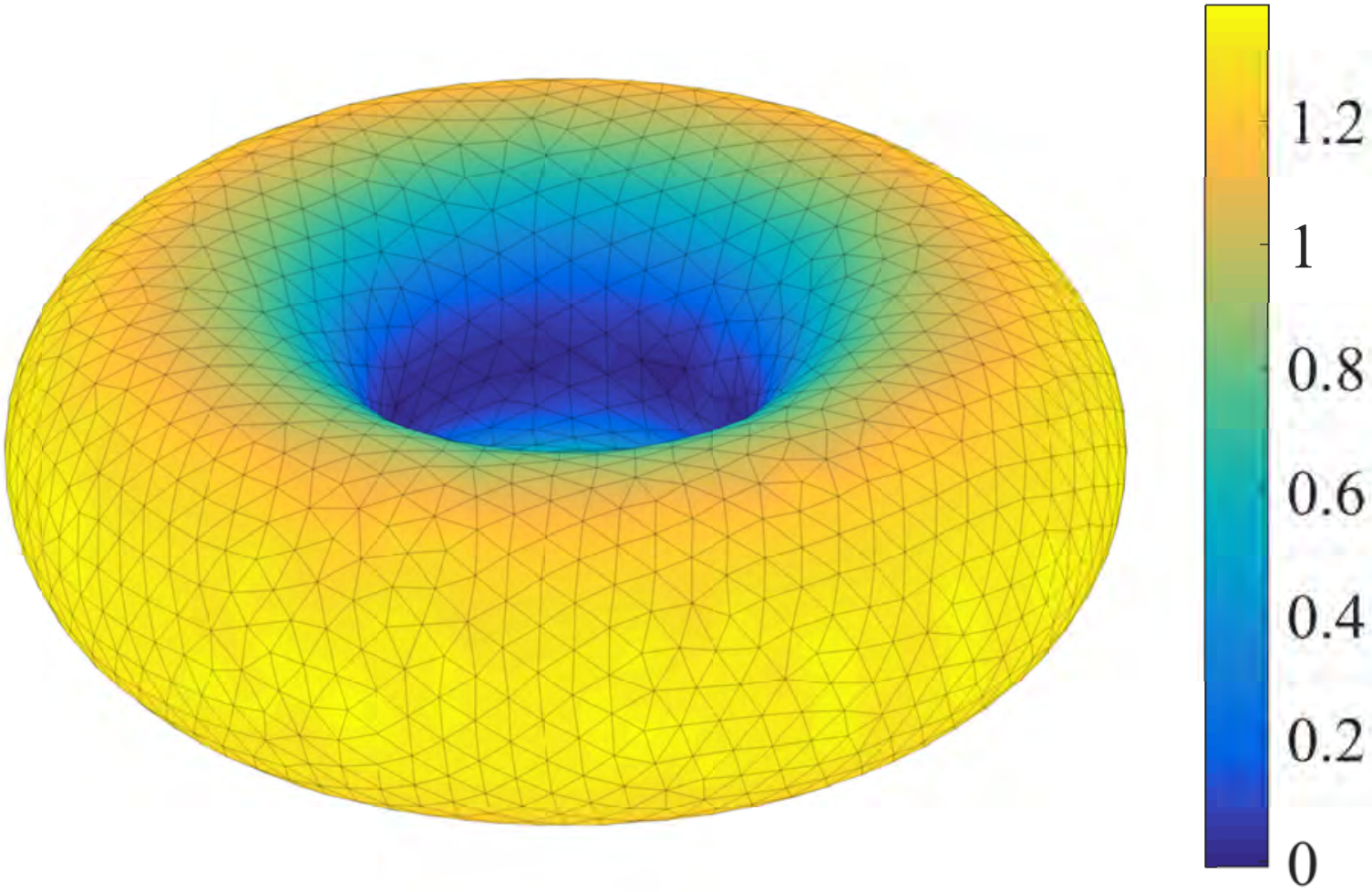}}
\par\end{centering}
\begin{centering}
\caption{Mean curvature on a torus with $a=1$, $h=0.0078$ (16384) vertices
and a unstructured mesh.\label{fig:MeanCurv}}
\par\end{centering}
\end{figure}

\begin{figure}
\centering{}\includegraphics[width=0.5\columnwidth]{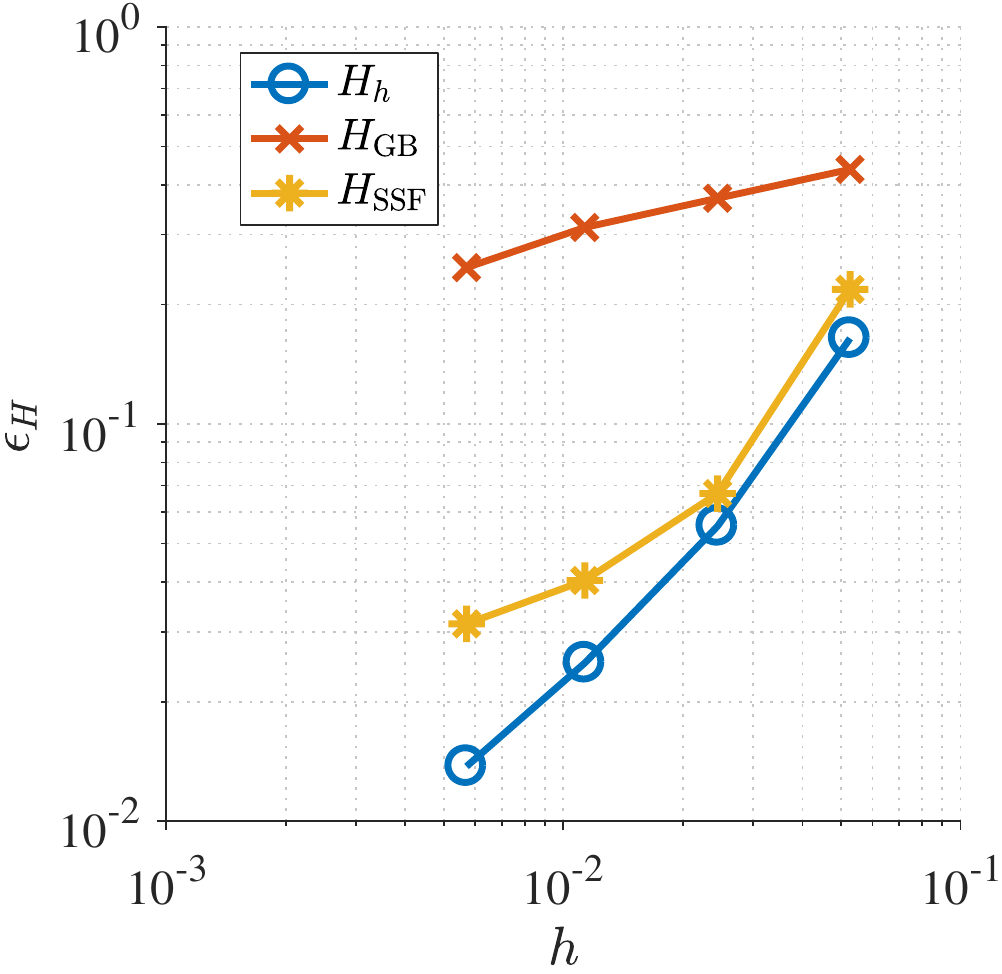}\caption{Mean curvature errors for a torus with $a=1$. $H_{h}$ computed with
$\gamma_{H}=0.05$ for all $h$.\label{fig:Mean-curvature-errors}}
\end{figure}

In Figure 
\ref{fig:MeanCurv} we give iso-plots of the mean curvature. Figure \ref{fig:Mean-curvature-errors}
shows the convergence of mean curvature.

\begin{table}
\begin{centering}
\begin{tabular}{|c|c|c|c|c|c|c|c|}
\hline 
$h$ & $\gamma_{H}$ & $H_{h}$ & $H_{h}$ rates & SSF & SSF rates & DLLB & DLLB rates\tabularnewline
\hline 
\hline 
0.0527 & 0.05 & 0.1643 & - & 0.2183 & - & 0.4377 & -\tabularnewline
\hline 
0.0245 & 0.05 & 0.0553 & 1.4178 & 0.0668 & 1.5412 & 0.3703 & 0.2179\tabularnewline
\hline 
0.0113 & 0.05 & 0.0250 & 1.0302 & 0.0404 & 0.6547 & 0.3126 & 0.2200\tabularnewline
\hline 
0.0057 & 0.05 & 0.0137 & 0.8746 & 0.0314 & 0.3670 & 0.2477 & 0.3394\tabularnewline
\hline 
\end{tabular}
\par\end{centering}
\centering{}\caption{Mean curvature errors for torus with $a=1$.}
\end{table}

\subsubsection{Stabilization sensitivity}

\begin{figure}
\begin{centering}
\includegraphics[width=1\columnwidth]{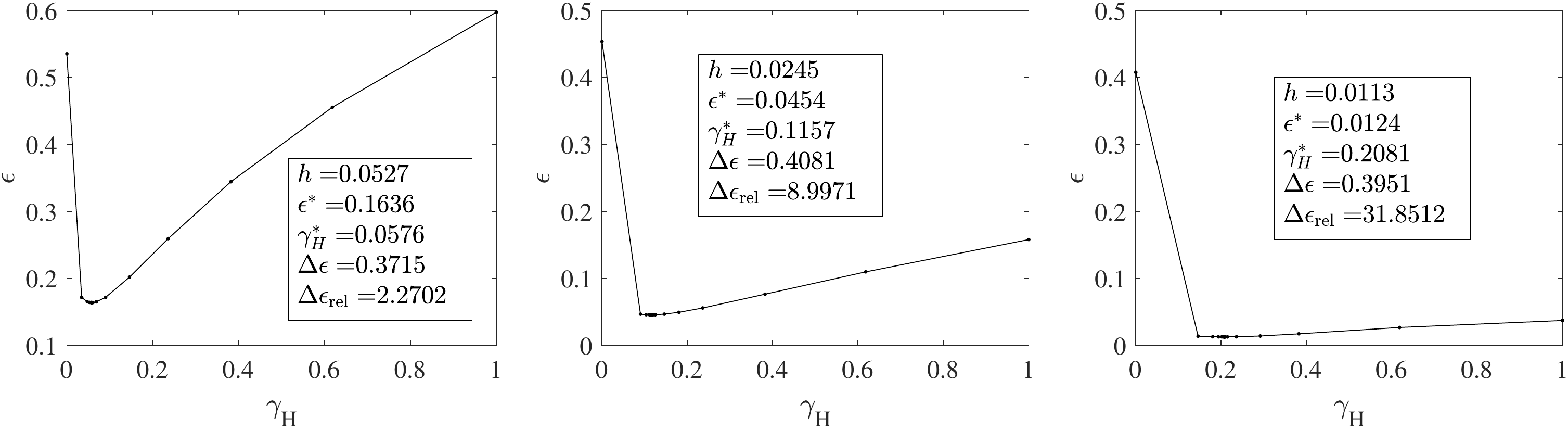}\caption{Mean curvature error $\epsilon_{H}(H_{h})$ as defined in (\ref{eq:MeanCurvError})
with respect to the stabilization factor $\gamma_{H}$.\label{fig:Mean-curvature-error-gammaH}}
\par\end{centering}
\end{figure}

We begin by analyzing how sensitive the mean curvature approximation
is with respect to the stabilization factor $\gamma_{H}$. The mean
curvature error is defined by

\begin{equation}
\epsilon_{H}=\|H_{\mathrm{exact}}-H_{h}\|_{L^2_h},\label{eq:MeanCurvError}
\end{equation}
where $H_{\mathrm{exact}}$ is the exact mean curvature computed on
a torus with $R=1,$ $r=\frac{1}{2}$ and $a=1$. A golden search
method is used to find the optimal stabilization factor and subsequently
to numerically analyze the impact of the stabilization choice with
regards to the error and mesh-size. 

\[
\begin{cases}
\underset{\gamma_{H}}{\mathrm{min}} & \epsilon_{H}\\
\mathrm{s.t.} & \gamma_{H}^{0}\leq\gamma_{H}\leq\gamma_{H}^{1}
\end{cases},
\]
where $\gamma_{H}^{0}=0$ and $\gamma_{H}^{1}=0.15$. The result of
this optimization is a discrete function of mean curvature error with
respect to the stabilization factor, $\epsilon_{H}(\gamma_{H})$,
see Figure \ref{fig:Mean-curvature-error-gammaH}. In Table \ref{tab:gammaOfH-Curvature}
we present the optimal stabilization factors and differences in mean
curvature error defined as $\Delta\epsilon_{H}=\epsilon_{H}(\gamma_{H}^{0})-\epsilon_{H}(\gamma_{H}^{*})$,
where $\gamma_{H}^{0}=0$ and $\gamma_{H}^{*}$ is the value of gamma
that minimizes $\epsilon_{H}$. Notice how the curves become more
planar, i.e., choosing a $\gamma_{H}$ that improves the solution
becomes less sensitive with the decrease in $h$.

\begin{table}
\begin{centering}
\begin{tabular}{|c|c|c|c|c|c|}
\hline 
$h$ & $\gamma^{*}$ & $\epsilon_{H}^{*}$ & $\Delta\epsilon_{H}$ & relative change & $\epsilon$ rate\tabularnewline
\hline 
\hline 
0.0527 & 0.0576 & 0.1636 & 0.3715 & 2.2702 & -\tabularnewline
\hline 
0.0245 & 0.1157 & 0.0454 & 0.4081 & 8.9971 & 1.6736\tabularnewline
\hline 
0.0113 & 0.2081 & 0.0124 & 0.3951 & 31.8512 & 1.6770\tabularnewline
\hline 
\end{tabular}
\par\end{centering}
\centering{}\caption{Stabilization factor $\gamma_{H}$ as a function of mesh-size $h$.\label{tab:gammaOfH-Curvature} }
\end{table}

\section*{Acknowledgements}

This research was supported in part by the Swedish Foundation for
Strategic Research Grant No. AM13-0029, the Swedish
Research Council Grants Nos. 2011-4992, 2013-4708,
and the Swedish Research Programme Essence.


\bibliographystyle{abbrv}
\bibliography{ComputGeom}
\newpage

\end{document}